\documentclass[11pt,leqno]{article}
\usepackage{amsthm,amsfonts,amssymb,epsfig,graphics,amsmath,oldgerm}
 \usepackage[latin1]{inputenc}\relax

\newcommand{\RR}{{\mathbb R}}

\newcommand{\ZZ}{{\mathbb Z}}
\newcommand{\TT}{{\mathbb T}}
\newcommand{\CC}{{\mathbb C}}
\newcommand{\EE }{{\mathbb E}}
\newcommand{\BB}{{\mathbb B}}

\newcommand{\FF}{{\mathbb F}}

\newcommand\cN{{\cal  N}}
\newcommand\cE{{\cal  E}}
\newcommand\cF{{\cal  F}}

\newcommand\cT{{\mathcal T}}
\newcommand\cS{{\mathcal S}}

\newcommand\vp{\varphi  }

\newcommand\D{\partial }
\newcommand\adots{\mathinner{\mkern2mu\raise1pt\hbox{.}
\mkern3mu\raise4pt\hbox{.}\mkern1mu\raise7pt\hbox{.}}}

\newcommand{\Id}{{\rm Id }}
\newcommand{\im}{{\rm Im }\, }

\newcommand{\la}{\langle }
\newcommand{\ra}{\rangle }

\newcommand{\rx}{{\mathrm{x } }}

\newcommand{\mez}{\frac{1}{2}}

\newcommand{\uE}{\underline{E}}

\newtheorem{theo}{Theorem}[section]
\newtheorem{prop}[theo]{Proposition}
\newtheorem{cor}[theo]{Corollary}
\newtheorem{lem}[theo]{Lemma}
\newtheorem{defi}[theo]{Definition}

\newtheorem{rem}[theo]{Remark}

\numberwithin{equation}{section}

\title{Space Propagation of Instabilities in  Zakharov Equations }

 \author{ Guy M\'etivier \footnote{Universit\'e Bordeaux 1, IMB, 
33405 Talence cedex, France; 
 Guy.Metivier@math.u-bordeaux1.fr} 
 \footnote{The first draft of this paper was elaborated during a visit at  the Ecole Polytechnique Fédérale de Lausanne, during the special semester devoted to 
 fluid dynamics. The authors thanks the organizers for their invitation and the EPFL  for its warm hospitality.  }}

\begin{document}

\maketitle

\begin{abstract}
In this paper we study an initial boundary value problem for 
Zakharov's equations, describing the space propagation of 
a laser beam entering in a plasma. We prove a strong instability 
result and prove  that the mathematical problem is ill-posed
in Sobolev spaces. We also show that it is well posed in spaces of  
analytic functions. Several consequences for the 
physical consistency of the model  are discussed.

\end{abstract} 

\tableofcontents


\section{Introduction}

Zakharov's
equations \cite{zakharov}  model   electronic plasma
waves,  describing  the coupling between  the slowly varying envelope of the electric 
field $E $  and  the low-frequency variation of
the density of the ions $ n$.  A commonly used form of the equations reads
\begin{equation}
\label{eq1}
\left\{ \begin{aligned}
& i( \D_t  + \frac{k_0 c^2}{\omega_0}\D_z)A + \frac{c^2}{2\omega_0}\Delta_x A =\frac{\omega_{pe}^2}{2n_0\omega_0} \delta n\,  A,\\
& (\D_t^2 - c_s^2\Delta_x) \delta n = \frac{\omega_{pe}^2}{4\pi m_i c^2}\Delta_x \vert A \vert^2,
\end{aligned}
\right.
\end{equation}
where $\omega_0$ is the frequency of the laser, $k_0$ its wave number and $\omega_{pe}$ the plasma electronic frequency; they  are linked by
the dispersion relation  $\omega_0^2 = \omega_{pe}^2+k_0^2c^2$ where   $c$ is the speed of light;  $n_0$ the mean density of the plasma, $m_i$ is the mass of the ions and $c_s$ the sound
velocity in the plasma. The space variables are  $(z, x)$, 
$z \in \RR$ and $x \in \RR^2$; $z$ is the direction of propagation of the laser beam
 and $x$ are the directions transversal  to the propagation. In this model, the transversal 
 dispersion is neglected. 

Reduced to  a dimensionless form (see in Section~6),  the equations become 
\begin{equation}
\label{1}
\left\{ \begin{aligned}
& i(\epsilon  \D_t  + \D_z)E + \Delta_x E = n E,\\
& (\D_t^2 - \Delta_x) n = \Delta_x \vert E \vert^2
\end{aligned}
\right.
\end{equation}
See \cite{Z2}
or \cite{sulem1} for  the introduction of  this kind of models for numerical simulation. 
  
The local in time Cauchy problem for 
\eqref{1} is now well understood: 
in \cite{LPS}, it is proved that it is  well posed, locally in time, for data in suitable  Sobolev spaces.  
This extends  previous results  for the classical Zakharov system, where  transversal dispersion 
is taken into account, that is when $\Delta_x$ is replaced by $\Delta_{(z, x)}$ (see \cite{sulem1, GTV, OT} and references therein). 
However, the system \eqref{1} is  quite different from the classical Zakharov system 
since the Cauchy problem for \emph{periodic}  data  exhibits strong instabilities 
of Hadamard's type (\cite{CM}).  Note that periodic data are considered in the quoted  paper 
as a model for data which do not vanish at infinity. 
We refer to  \cite{CM}
for a discussion of the physical relevance of the different frameworks.  
In particular, we  recall  in Section~6,  why the periodic context is well adpated to 
 physical situations, where the envelop of the beam has rapid oscillations (speckles).

In this paper we consider   a boundary value problem for \eqref{1} which  
 models the the propagation of 
a laser beam entering the plasma at an interface $\{ z = 0 \}$. This approach is very common 
in physics where people is actually more interested in describing \emph{ the propagation 
in space} rather than in time, i.e.  considering $z$ as the propagation variable. Indeed, this  
is an  underlying idea in many of the paraxial approximations, like the one which yields 
the Schr\"odinger equation in \eqref{1}. This idea is also very common 
in numerical simulations.   In this case, the system \eqref{1} is considered 
in the half space $\{ z \ge 0 \}$, for positive times $\{ t \ge 0 \}$ together with 
initial-boundary conditions 
\begin{equation}
\label{1.2}
n _{| t = 0}  = 0, \quad \D_t n_{| t = 0} = 0, \quad  E_{ z = 0}  = E_0 , 
\quad E_{t = 0} = 0 . 
\end{equation}
This approach is also very natural when 
the parameter $\epsilon$ is small.
The solutions are expected to vanish for $t \le z$, by finite speed of propagation. 
Changing $t$ to $t - \epsilon z$, the system reads; 
\begin{equation}
\label{1.3}
\left\{ \begin{aligned}
& i  \D_z E + \Delta_x E = n E,\\
& (\D_t^2 - \Delta_x ) n = \Delta_x \vert E \vert^2, 
\end{aligned}
\right.
\end{equation}
 \begin{equation}
\label{1.4}
\left\{\begin{aligned}
& n _{| t = 0}  = 0, \quad \D_t n_{| t = 0} = 0, 
\\
& E_{ z = 0}  = E_0 .
\end{aligned}\right. 
\end{equation}
We look for solutions $U = (E, n)$ of \eqref{1.3} \eqref{1.4}  which are periodic in $x$ 
with period $L$. We denote by $\TT $ the 
corresponding torus  $(\RR/ 2 \pi L)^2$. 
As recalled in Section 6,  periodicity is somewhat  natural when the envelop of the laser beam has 
a transversal structure  which involves length-scales that are large with respect to the 
length-scale of the laser but small compared to the width of the beam (see \cite{CM}). 
Periodicity is also natural for the use of  spectral methods in numerical simulations. 

The main result of this paper is a strong instability result for \eqref{1.3}, in the spirit of \cite{CM}. 
It says  that  \emph{the boundary value problem \eqref{1.3} \eqref{1.4} is ill-posed in Sobolev
spaces}.  Prior to this result, we make a detailed analysis of \emph{amplification properties} of the linearized   equations : at a  space frequency of length $k$, the amplification is exponential in 
$e^{ c \sqrt k}$. This implies  a Hadamard's type instability and the ill-posedness in Sobolev 
spaces, as well a the local well-posedness in spaces of real analytic functions. 
For physical interpretations the frequencies must be considered as  bounded
while  for numerical applications they are filtered by one way or another.
Thus is is important to give a numerical value to  the amplification rate  and to be 
evaluate the amplification effect on definite  values of the fields. 
This is done in Section~6.

\bigbreak 

 Let us describe now the main contents of the paper.  We first consider the linearized equations.
 Consider here a  constant solution 
\begin{equation}
\label{1.6}
\underline U =  (  \underline E,  0) , \quad   \underline E \ne 0, 
\end{equation} 
which satisfies \eqref{1.3} \eqref{1.4} with boundary data $E_0 = \underline E$. 
The homogeneous linearized equations near $\underline U$,  
at the frequency $\xi \in \RR^2$ read:   
 \begin{equation}
\label{1.6b}
\left\{ \begin{aligned}
& i \D_z  u   - k^2    u  -    \underline E   n =  0   ,\\
& i \D_z  v   +  k^2    v   +    \overline{ \underline E }  n =   0   , \\ 
& (\D_t^2 + k^2 )  n + k^2   \underline  E   ( u  + v  ) =  0 , 
\end{aligned}
\right.
\end{equation}
where  $k := | \xi |$ and $u$ [resp. $v$] [resp $n$] denote the 
the Fourier coefficient at the frequency $\xi$ of the variation of the field $ E $ [resp. $\overline E$]
[resp of the density]. They are  supplemented  with  initial boundary boundary conditions:   
 \begin{equation}
\label{1.7b}
n _{| t = 0}  = 0, \quad \D_t n_{| t = 0} = 0, \quad  u_{ z = 0}  = u_0, \quad 
 v_{ z = 0}  = v_0 .
\end{equation}
The equations \eqref{1.6b} \eqref{1.7b} form a well posed hyperbolic Goursat problem 
(see \cite{Al}).  Our concern is to understand its behavior for large $k$. 
The fundamental solution  is studied in details in Section~2, where we also construct blowing up solutions: 
   
\begin{theo}
\label{thampli}
 There are initial data $ | u_0 | \le 1 $ and $v_0 = 0$, such that for $k $ large, 
 $t \ge \mez$ and $\rho :=   \sqrt{ 2k |\uE|^2 zt } $ large,  the solution of the 
 homogeneous equation \eqref{1.6b} satisfies  
 \begin{equation}
 | n(t, z) | \gtrsim  k  \rho^{- \frac{5}{2}} e^{ \rho} . 
 \end{equation} 
\end{theo}

\medbreak
 
The exponential amplification in $e^{ c \sqrt k} $ is the signal of a strong instability.  
We   construct solutions on domains 
\begin{equation}
\label{1.7}
\Omega = \{  (t, z, x ) \in [0, T] \times [0, Z] \times \TT^2  \ ; \   z t \le  \delta \} 
\end{equation}
and prove that they  blow up on the part of the  boundary 
\begin{equation}
\label{1.8}
\Gamma = \{  (t, z, x) \in \Omega  \ ; \   z t  =   \delta \} 
\end{equation}
which is not empty if $TZ  >   \delta$. 

\begin{theo}
\label{mainth}
For all  $s$, $T > 0$ $L > 0$  and $\underline E \ne 0$, there are sequences 
$\delta_k$ and $Z_k $ and families of  solutions  $U_k  = \underline U +  (e_k, n_k )$ of \eqref{1.3}, in $C^0( \Omega_k; H^s(\TT))$
such that 
\begin{eqnarray}
&& \|  e_k{}_{| z = 0 }    \|_{H^s([0, T ] \times  \TT^2)}  \ \to \ 0,  
\\
&& Z_k \to 0,  \quad \delta_k \to 0, \quad  TZ_k  >   \delta_k, 
\\
&&    \sup_{(t,z) \in \Gamma_k}   \|    n_k( t,z)  \|_{L^2(\TT)}  \ \to \  \infty . 
\end{eqnarray}

\end{theo}

This theorem is proved in Section~ 3, with technical details postponed to Section~4. 
This nonlinear instability result is pretty strong:  not only the 
 {amplification}    $ \|  e   \|_{L^2}  / \| e_{| z = 0}  \|_{H^s}  $ is arbitrarily large, in 
arbitrarily small distance  $Z$, with arbitrary loss of derivatives $ s$, but there is 
an effective \emph{blow up} of the   norm of $n_k$.

This  analysis reveals the importance of the amplification factor, 
 $\rho :=   \sqrt{ 2k |\uE|^2 zt } $, and indeed there is a good uniform stability 
 for a filtered system at frequencies $| \xi | \le k$, 
 on the domain $  \{  (t, z, x ) \in [0, T] \times [0, Z] \times \TT^2  \}$  provided that 
 \begin{equation}
 \label{a1.15}
    k Z T | \uE |^2  \ll    1 . 
 \end{equation}

\medbreak
Instead of filtering the frequencies, another mathematical approach is to counterbalance  the amplification  by an exponential decay
of the Fourier coefficients.  This means that one works in spaces of real analytic functions. 
In this framework,  we   prove in Section~5  a local existence theorem by an easy  adaptation of 
proof of the existence of analytic solutions to the hyperbolic Goursat problem (see \cite{Wag}). 
 The interesting point is that the length $Z$ of propagation satisfies  
 an estimate which is very similar to \eqref{a1.15}: 
  \begin{equation}
  \label{a1.16}
    R  Z T \| \uE_0 \|_{\EE_R} ^2  \ll    1 . 
 \end{equation} 
where $\EE_R$ is a space of analytic functions for the boundary data on 
$[0, T] \times \TT^2$ and $R^{-1}$ measures the width of the complex domain 
where $E_0$ can be extended. 
Indeed, for $E_0 (t,x) = e_0(t) e^{ i \xi \cdot x}$ with $| \xi | = k$, 
 and  $R \approx k$,   there holds 
 $\| E_0 \|_{\EE_R}  \lesssim   \uE   =  \| e \|_{L^\infty} $
 and \eqref{a1.16} is equivalent to \eqref{a1.16}. 
  
 \emph{ This shows that \eqref{a1.15} or \eqref{a1.16} can be seen as stability criteria for the 
  Goursat problem for Zakharov equations \eqref{1}.}
  
\medbreak

Section~6 is devoted to a  qualitative discussion of the results.  We discuss several points.

\quad -  Taking  $\underline E$ to be a constant is not physically realistic : it would mean that 
the envelop of the mean electric field has a jump.    The case where   $\underline E$ is a 
smooth function of  time  will be  
briefly  discussed. It yields additional technical difficulties but does not change qualitatively the results. 

\quad -  The boundary data in Theorem~\ref{thampli} are very particular. Thus, it is important to  understand better how general solutions behave. Actually, the rate of amplification 
also depends on the time frequencies. The underlying phenomenon is a \emph{resonance}  between    plasma waves governed by the  Schr\"odinger equation and  acoustic waves for $n$. 
In other words, \emph{only the acoustic oscillations}  $e^{ -i \omega t + i \xi x}$
with $\omega^2 = | \xi|^2$  which present in the boundary data $E_0$  are 
effectively amplified at the given exponential rate. 
 Mathematically, in general, these acoustic frequencies have a nonvanishing amplitude  because the 
 signals exactly vanish in the past and thus their Laplace-Fourier transform 
 has no lacuna. Physically, the  acoustic frequencies,  even when   absent from  the main scene, 
 can be present in  a background noise. 
 By a standard plane wave analysis, we will also give an amplification rate 
 for  oscillations  which are not exactly acoustic, but this  is not a correct approach for the 
 Gourset problem. 
 
\quad - The blow up in Theorem~\ref{mainth} is totally unphysical, since $n$ is a variation of density and thus must
remain bounded. 
However, for physical interpretations,  one must keep in mind that 
that \eqref{eq1} is only a \emph{model}  which has a limited range of validity.
In particular,  it is tacitly assumed that the variations of the ion density are not 
too large and that  the paraxial approximation for the envelop  is valid. 
Moreover, the   physical   frequencies $k$ are   \emph{bounded}.   
One crucial question 
is to know  \emph{  wether the  factor $\rho =  \sqrt{2 k | \uE |^2 z t}$ is   small or large}. 
In Section~6, we will give standard physical data for laser-plasma propagation
showing that this \emph{factor can be  large } ($10$ to $10^2$) for not very intense fields 
$E $ of order $10^9$ or $10^{10}  \ \mathrm{W m^{-1}}$. In this case,  that 
the amplification $e^  \rho $ ranges from $10^4 $ to $10^8$. Thus, the acoustic boundary  oscillations 
can be ignored only if their relative value is $\ll e^{ - \rho}$, which is much beyond
 the usual admissible errors.  This seems to indicate that the model, as is it, 
 is not well adapted to the propagation of intense laser beams. 
 
\quad - For numerical simulations, the analysis shows that increasing  the number of Fourier modes, which is natural 
to improve the accuracy of computations, may introduce  strong instabilities.  Of course, one can eliminate most of them by filtering out the 
bad acoustic oscillations, but then the question is the relevance of the computations with 
respect to the model.



\section{The linear instability}

Consider the linearized equation from \eqref{1.3} around $(\underline E,   0 )$ : 
\begin{equation}
\label{2.1}
\left\{ \begin{aligned}
& i \D_z  e  + \Delta_x   e -    \underline E   n = f  ,\\
& (\D_t^2 - \Delta_x)  n -  \Delta_x (   \overline {\underline E}  e + \underline E \overline e ) =  \Delta_x h  
\end{aligned}
\right.
\end{equation}
\begin{equation}
\label{2.2}
n _{| t = 0}  = 0, \quad \D_t n_{| t = 0} = 0, \quad  e_{ z = 0}  = e_0 . 
\end{equation}
 Multiplying $e$ by a constant phase factor 
 $e^{ i \theta}$, there is no restriction  in assuming that $\underline E$ is real. Taking 
 $u = e$ and $v = \overline e$ as (independent) unknowns, the system reads
 \begin{equation}
\label{2.3}
\left\{ \begin{aligned}
& i \D_z  u  + \Delta_x   u -    \underline E   n = f  ,\\
& i \D_z  v  -  \Delta_x   v  +      \underline E     n =    g   , \\ 
& (\D_t^2 - \Delta_x)  n -  \underline E   \Delta_x  ( u + v  ) = \Delta_x  h , 
\end{aligned}
\right.
\end{equation}
with $ g = - \overline f$. 

 Performing  a Fourier series expansion in $x$ (or a Fourier transform), 
 and still denoting by $u$, $v$ and $n$ the Fourier coefficients, the 
 equations at the frequency $\xi \in \RR^2$ read  with $k := | \xi |$:  
 \begin{equation}
\label{2.4}
\left\{ \begin{aligned}
& i \D_z  u   - k^2    u  -    \underline E   n =  f   ,\\
& i \D_z  v   +  k^2    v   +    \overline{ \underline E }  n =   g   , \\ 
& (\D_t^2 + k^2 )  n + k^2   \underline  E   ( u  + v  ) = - k^2  h , 
\end{aligned}
\right.
\end{equation}
together with  initial boundary boundary conditions:   
 \begin{equation}
\label{2.5}
n _{| t = 0}  = 0, \quad \D_t n_{| t = 0} = 0, \quad  u_{ z = 0}  = u_0, \quad 
 v_{ z = 0}  = v_0 .
\end{equation}
For $L$-periodic functions, the frequencies $\xi \in \frac{2\pi}{L} \ZZ^2$. 

The Goursat problem 
 \eqref{2.4} \eqref{2.5}  is well posed (\cite{Al}). 
   The main purpose of  this section is to prove estimates for the fundamental solutions
(Propositions~\ref{prop2.5} and \ref{prop2.6} below), and give an example
of a solution of the homogeneous equation which is amplified at   indicated rate indicated 
in  Theorem~\ref{thampli}. 
(see also Theorem~\ref{th2.7} below for a more precise statement).


\subsection{The fundamental solution} 

 Extend the functions by $0$ for $t < 0$ and 
perform a  Fourier-Laplace transform in time; this amounts to replace 
 $\D_t$ by $   i \zeta$  with $\zeta$ lying in the lower half plane 
 $\{  \im \zeta < 0 \} $:   the third equation in \eqref{2.4} and the homogeneous initial conditions for $n$ imply  that   
 \begin{equation}
 \label{2.6}
 \hat n   =    \frac{ k^2} { \zeta^2  -  k^2}   \big(  \uE ( \hat u +  \hat v  )   + \hat h \big) . 
 \end{equation}  
 We denote here by $\hat \phi ( z, \zeta)$ the Fourier-Laplace transform of $\phi(z, t)$. 
 We end up with the system
 \begin{equation}
 \label{2.7}
    \D_z \hat U  =   i   A \hat  U  +  \hat F , \qquad  \hat U_{| z = 0} = \hat U_0
 \end{equation}
 for 
 \begin{equation}
 \label{2.8}
 \hat U = \begin{pmatrix}
    \hat u     \\
  \hat v    
\end{pmatrix}
,  \qquad  \hat F  = \begin{pmatrix}
    \hat f     \\
\hat g     
\end{pmatrix}   +  i   \frac{ k^2  \uE \hat h } { \zeta^2  -  k^2}
\begin{pmatrix}
       -1      \\
  1     
\end{pmatrix} , 
 \end{equation}
 with 
 \begin{equation}
 \label{2.9}
A    (\zeta) =    \begin{pmatrix}
     -  k^2 -   a  &   -   a    \\
    a    &     a +   k^2
\end{pmatrix}, \qquad  a = \frac{ k^2 | \underline E |^2}{ \zeta^2 - k^2}.  
 \end{equation}

 Note that $A(\zeta)$ is bounded and holomorphic for $\{ \im \zeta \le - \gamma \}$ for 
 all $\gamma > 0$ and therefore, by inverse Laplace transform: 
 \begin{lem}
 \label{lem2.1}
 The solution of the homogeneous system $\eqref{2.4}$ with initial-boundary values
\eqref{2.5} where  $U_0 = {}^t (u_0, v_0)  \in C^\infty_0 (\overline \RR_+) $ on $\{ z = 0 \}$, is 
 \begin{eqnarray*} 
&& U (t, z) 
 =   \frac{1}{2 \pi} 
\int_{\RR - i \gamma}  e^{  i \big( t \zeta + z A  (\zeta)  \big)  }  \hat U_0 (\zeta) \  d \zeta , 
\\
&&n (t, z) 
 =   \frac{1}{2 \pi} 
\int_{\RR - i \gamma}  \ell \cdot  e^{  i \big( t \zeta + z A  (\zeta)  \big)  }  \hat U_0 (\zeta)
 \frac{ k^2  E    } { \zeta^2  -  k^2}   d \zeta , 
 \end{eqnarray*}
   where $\ell := (1, 1 )$ and $\gamma $ is any positive real number.
 \end{lem}
   
   Note that the integrals above are convergent, as a consequence of the estimates given 
   below. 
   
   \medbreak

 The fundamental solution is therefore linked to the distribution 
  \begin{equation}
 \label{2.10} 
\mathcal{E} (t, z) 
 =   \frac{1}{2 \pi} 
\int_{\RR - i \gamma}  e^{  i \big( t \zeta + z A  (\zeta)  \big)  }\  d \zeta , 
\end{equation}
where the integral is taken in the sense of  an inverse Laplace transform.
More generally, we are led to consider integrals of the form 
   \begin{equation}
 \label{2.11} 
\mathcal{I } (t, z) 
 =   \frac{1}{2 \pi} 
\int_{\RR - i \gamma}  e^{  i \big( t \zeta + z A  (\zeta)  \big)  }   p(\zeta) \  d \zeta , 
\end{equation}
for   rational functions $p$. 

 Note that $a$ and $A$  and thus 
  $e^{  i \big( t \zeta + z A (\zeta)  \big) } $ are  holomorphic 
  in $\CC \setminus\{ -k, +k\}$.  In addition, $a   (\zeta ) = O( | \zeta |^{-2})$ 
  at infinity,  and thus  
  \begin{equation}
  \label{2.11n}
  A(\zeta ) - A(\infty)  = O( | \zeta |^{-2}) , \qquad 
  e^{ i z A(\zeta )}  - e^{ i z A(\infty)}   = O( | \zeta |^{-2})
  \end{equation}
  where 
  \begin{equation}
  \label{2.12}
  A( \infty) =  \begin{pmatrix}
     -  k^2   &   0     \\
    0     &       k^2
\end{pmatrix}. 
  \end{equation}
  
  This implies the following

  \begin{lem}
  \label{lem2.2}
       If   $p $  is a rational function with poles in $\{ \im  \zeta \ge 0 \}$,   and 
       $ p(\zeta)  = O(|\zeta |^{-2}) $  at infinity, 
       then the integral  in $\eqref{2.8}$ is absolutely convergent for all 
       $t \in \RR$. It vanishes for $t < 0$ and for $t \ge 0$,   
    \begin{equation}
    \label{2.13}
\mathcal{I}  (t, z)    =  \frac{1}{2\pi} \int_{\Gamma }  
e^{   i  t \zeta  + i z A   (\zeta )   }
 p  (\zeta)  \ d \zeta 
  \end{equation}
 where $\Gamma$ is any simple contour oriented positively
 winding around the poles of $p $ and $a$. 
  \end{lem}
  
  \begin{proof}
  The integral is clearly convergent and letting $\gamma $ tend to 
  $+ \infty$ implies that it vanishes when $t < 0$. 
  When $t \ge 0$, the integral over a large half circle in the upper half space 
  tends to $0$ as the radius tends to infinity, allowing to close the integration path. 
  \end{proof}

  Similarly,  we can split $\cE$ into two parts: 
  $$
  \mathcal{E} (t, z) 
 =   \frac{1}{2 \pi} 
\int_{\RR - i \gamma}  e^{  i \big( t \zeta + z A  (\infty)  \big)  }\  d \zeta 
+
 \frac{1}{2 \pi} 
\int_{\RR - i \gamma}  e^{  i   t \zeta } 
\big( e^{ i z A  (\zeta)} - e^{ i z A  (\infty)}   \big)  \  d \zeta . 
    $$
 Thanks to \eqref{2.11n}, the second integral can be deformed to an  integral on a closed contour $\Gamma$, 
 on which the integral of the entire function 
 $e^{  i   t \zeta }   e^{ i z A  (\infty)} $   vanishes. Therefore, 
    
     \begin{lem}
     \label{lem2.3}
       The distribution $\mathcal{E}$ defined in $\eqref{2.10}$  is equal to   
    \begin{equation}
    \label{2.14}
\mathcal{E}  (t, z)    = 
e^{ i z A(\infty)}  \ \delta_{ t = 0}   +  \cE_0(t, z) 
  \end{equation}
  where $\cE_0$  vanishes for $t < 0$  and is equal to 
  \begin{equation}
  \label{2.15}
  \cE_0(t, z) =  
 \frac{1}{2\pi} \int_{\Gamma }  
e^{   i  t \zeta  + i z A_k (\zeta )   }
 \ d \zeta 
\end{equation} 
when $ t \ge 0$,  where $\Gamma$ is any simple contour oriented positively
surrounding   the poles of  $a$. 
  \end{lem}


\subsection{Estimates} 
    We give sharp  upper bounds for  the contour integrals in \eqref{2.13} and \eqref{2.15}. 
     The matrix $A$ is traceless, therefore  $ A^2 =  - (\det A) \Id =  k^2 (k^2 + 2 a) \Id$
     and 
     $$
      e^{ i z A}   =  \cos (k z \lambda) \Id \   +\   \frac{\sin (k z \lambda)} { k \lambda}   A  
     $$ 
     where 
    \begin{equation}
    \label{2.16}
    \lambda = \sqrt { k^2 + 2 a}. 
    \end{equation}
 The choice of the square root is irrelevant in the expression above of $e^{ i z A}$.  
  Therefore, 
     \begin{equation}
\label{2.17} 
\begin{aligned}
    e^{ i t \zeta +  i z  A }   =   & 
    \frac{1}{4} e^{  i t \zeta  -   i k \lambda z} \begin{pmatrix}
    2 + \frac{\lambda}{k} + \frac{k}{\lambda}    &    \frac{\lambda}{k} - \frac{k}{\lambda}     \\
  \frac{k}{\lambda}  -  \frac{\lambda}{k}      &  2 -  \frac{\lambda}{k} -  \frac{k}{\lambda} 
\end{pmatrix}
    \\
    &   +    \frac{1}{4} e^{i t \zeta +   i k \lambda z} \begin{pmatrix}
    2 - \frac{\lambda}{k} - \frac{k}{\lambda}    &  \frac{k}{\lambda}  -  \frac{\lambda}{k}   \\
   \frac{\lambda}{k} - \frac{k}{\lambda}     &  2 + \frac{\lambda}{k} + \frac{k}{\lambda} 
\end{pmatrix}. 
       \end{aligned}
   \end{equation}    
   This reveals the phases 
   \begin{equation}
   \label{2.18}
   \vp = t \zeta  - z k  \lambda, \qquad 
 \psi = t \zeta  +  z k  \lambda
 \end{equation}

 \begin{lem}
 \label{lem2.4f}
For  $k \ge \max\{ 1,  4 | \uE |^2 \} $,    
  choosing the principal 
determination of the square root in $ \eqref{2.16}$ defines 
$\lambda(k + w) $ as an holomorphic function of  $w$ 
for  $1 \le | w |  \le \mez k$. Moreover, for   all $t \ge 0$, $z \ge 0$ and $w$ in this annulus, there holds 
\begin{eqnarray*}
&\Big|  \vp ( k + w )   -  \big(   k t  - k^2 z    + t w - {\displaystyle\frac{ k z |\uE|^2}{ 2w} } \big) \Big| \le  
   z \big( \frac{1}{3}  |\uE|^2  +  \frac{1}{5}  |\uE|^4 \big)   ,
\\
&\Big|  \psi  ( k + w )   - \big(  k t  +  k^2 z   + t w + {\displaystyle  \frac{ k z |\uE|^2}{ 2w}  } 
\big) \Big|  \le   z \big( \frac{1}{3}  | \uE |^2 +    \frac{1}{5}  |\uE|^4 \big). 
\end{eqnarray*}

\end{lem}

\begin{proof}
For $ | w |  \le \mez k$, 
$$
a(k + w)    = \frac{ k^2 |\uE|^2 } { w ( w + 2k) }    =  
 \frac{ k | \uE |^2}{2  w }      -     \frac{| \uE |^2}{ 4 (  1+  \frac{ w}{2k})  } 
$$ 
Thus 
$$
\Big |  a   -     \frac{  k  | \uE |^2}{2  w }   \Big |   \le   \frac{1}{3}  \  | \uE |^2 
$$
If in addition $   | w  | \ge 1 $ and $k \ge \max\{ 1,  4 | \uE |^2 \} $, we have 
$$
\Big| \frac{a}{k}   \Big| \le | \uE |^2 , \qquad   \Big| \frac{2 a}{k^2 } \Big| 
\le \frac{ 2 | \uE |^2} { k} \le \mez. 
$$
  Thus  we can  choose $\lambda$ to be 
the principal determination of the square root of $ k^2 + 2 a$ and 
\begin{equation}
\label{2.19}
k \lambda =  k^2  + a  + \frac{a^2}{k^2}  G \big( \frac{2a}{k^2} \big)
\end{equation}
where $G$ is holomorphic on the unit disc.  Substituting, we find that 
\begin{equation}
\label{2.20}
\big|  k \lambda -  k ^2 -   \frac{ k | \uE |^2}{ 2 w } \big|   \le   \frac{1}{3}  | \uE|^ 2 + 
| \uE |^4  \sup_{|\zeta  | \le \mez} | G ( \zeta ) |  .  
\end{equation}
The supremum is $| G (- \mez)| \le \frac{1}{5}$ and the lemma follows. 
\end{proof}

\begin{prop}
\label{prop2.5}

Then there is a constant $C$,  such that for all $\uE$, $k$, $t$ and $z$ satisfying
\begin{equation}
\label{2.22f}
 | \uE | \le 1 ,  \quad  k \ge \max\{ 1,  4 | \uE |^2 \}, \quad  
  t \in [0, 1], \quad  z \in [0, \mez], 
  \end{equation}
   the contour integral in $\eqref{2.15}$ 
  satisfies :  
\begin{equation*}
\big| \mathcal{E}_0 (t, z ) \big|  \le C \frac{ r   \ e^\rho}{ \sqrt {1+ \rho}} 
\end{equation*}
with $\rho = \sqrt { 2  k |E|^2 zt }$ and $r =  1 + \frac{ k | \uE |^2 z}{ 1 + \rho}$. 
\end{prop}
\medbreak
\begin{proof}
{\bf a) }    We choose  $\Gamma $ to be the union of the circles 
$$
\Gamma_1 = \{  | \zeta - k | = r \} , \qquad  
\Gamma_2  = \{  | \zeta + k | = r \} ,
$$
with  $r \in [1, \mez k] $ to be chosen depending on $ t$ and $z$. 

By  \eqref{2.20} and \eqref{2.22f}, we have the following bounds: 
 \begin{equation*}
\Big|  \frac{\lambda}{k} - 1 \Big| \le \frac{| \uE |^2}{2 k r} + \frac{ | \uE |^2 z}{ k^2} 
\le  \frac{| \uE |^2}{  k  }
\le \frac{1}{4} . 
 \end{equation*}
Furthermore
 \begin{equation*}
\Big|  \frac{k}{\lambda} - 1 \Big| \le \frac{| \uE |^2}{2 | \lambda|  r} + \frac{ | \uE |^2 z}{ k | \lambda| } 
\le   \frac{ 2 | \uE |^2}{  k  }    \le \mez . 
 \end{equation*}
 Thus 
\begin{equation*}
\Big|  \frac{\lambda}{k} + \frac{k}{\lambda}  \Big| \le 3, \quad
 \Big|  \frac{\lambda}{k} -  \frac{k}{\lambda}  \Big| \le   \frac{ 3  | \uE |^2 }{ k } , 
 \quad
\Big| 2 -   \frac{\lambda}{k} -  \frac{k}{\lambda}  \Big| \le   \frac{ 3  | \uE |^2 }{ k } . 
 \end{equation*}
 Therefore, all the entries of the two matrices present in \eqref{2.17} are bounded by $5$. 
 
 \medbreak
 {\bf b) }   If   
 \begin{equation}
 \label{2.21}
 k  | E |^2 z t \ge 2, 
 \end{equation}
 then we choose  
\begin{equation}
\label{2.22} 
  r = \sqrt{ \frac{k z |E|^2}{2 t} }. 
\end{equation}
 By \eqref{2.21} and \eqref{2.22f},  $r  \in [1, \mez k] $.  
 For 
 $\zeta = k + r e^{ i \theta} \in \Gamma_1$, Lemma~\ref{lem2.4f} implies that  the imaginary parts of the phases satisfy 
 \begin{eqnarray}
 \label{2.23}
&& \big|   \im  \vp ( k + r e^{i \theta}  )   - 
 \rho \sin \theta  \big|  \le   z  | \uE |^2  ,      
\\
&&
\label{2.24} 
\big|   \im  \psi  ( k + r e^{i \theta}   \big| \le  z  | \uE |^2  . 
\end{eqnarray}

  Therefore,  the integral over $\Gamma_1$ contributes to $\cE_0$ 
  to a matrix whose entries are bounded by 
   $$
 \frac{   1}{2 \pi}    \int_0^{2 \pi}  
 \frac{5}{4} e^{ z | \uE |^2}  \big( 1+    e^{ - \rho \sin \theta}   \big) \ r  \  d \theta. 
 $$
By  symmetry, the integral over $\Gamma_2$ is estimated similarly, and therefore, 
for $\rho \ge 2$, 
the entries of $\cE_0$ are bounded by 
  $$
 \frac{ 5   r }{4  \pi}    e^{ z | \uE |^2}  \int_0^{2 \pi}  
   \big( 1+    e^{ - \rho \sin \theta}   \big)    \  d \theta. 
 $$
 The integral is bounded by   
 $$
 \begin{aligned}
 3 \pi + \int_{- \pi/2}^{\pi/2} e^{ \rho \cos \theta} \ d \theta & = 
 3 \pi +   2  e^{ \rho} \int_{ - \pi/4} ^{ \pi/4} e^{ -  2 \rho \sin^2 y } \ dy 
 \\
 & \le 3 \pi +   2  e^{ \rho} \int_{ - \sqrt 2 /2 } ^{ \sqrt 2/2 } e^{ -  2 \rho t^2  } \ \frac{dt}{\sqrt{1 - t^2} } 
 \le 3 \pi + \frac{2 \sqrt \pi e^{ \rho}}{ \sqrt \rho} . 
 \end{aligned}
 $$
 Therefore, for $\rho \ge 2$,  the entries of $\cE_0$ are bounded by 
  \begin{equation}
\big| \mathcal{E}_0 (t, z ) \big|  \le  5 \ e^{ z | \uE |^2} \   \sqrt { \frac{ | \uE |^2 kz}{t}} \ 
\Big(  \frac{3}{4}  +   \frac{    e^\rho }{2 \sqrt{ \pi \rho} }\Big) 
\le     5 \ \sqrt e   \sqrt { \frac{ | \uE |^2 kz}{t}} \ 
    \frac{    e^\rho }{  \sqrt{ \pi \rho} } 
\end{equation}

     \medbreak
     {\bf c) }   If 
      \begin{equation}
      \label{2.25}
\mez \rho^2 =   k  | E |^2 z t \le 2, 
 \end{equation}
  we choose 
  \begin{equation}
  \label{2.26}
  r = \max\{ 1 , k | E |^2 z\}    \in [1, \mez k] .
  \end{equation}
     Consider first the integral on $\Gamma_1$.  The entries of the corresponding  matrix 
     The imaginary  part of the phases are bounded by
       \begin{equation*}
       \big| \im \vp \big| \le  t r  + z  k \frac{| \uE |^2}{r} \le \max\{ t  ,  t k | E |^2 z\}  + 1 
   \le 2 + \rho. 
       \end{equation*}
 Therefore, the  integral are bounded by 
 $$
 \frac{ 5 r}{4} e^{ 2 + \rho} \le 10 \ e^\rho. 
 $$
The analysis of the integral over $\Gamma_2$ is similar, 
     and combining the estimates above, one obtains \eqref{2.18}
     \end{proof}

     \begin{prop}
     \label{prop2.6}
Denote by $\mathcal{E}_1$ [resp.  $\mathcal{E}_2$ ]   the contour integral in $\eqref{2.11}$ with 
$ p  = \frac{k^2}{ \zeta^2 - k^2}$ [resp. $ p  = \frac{k^4}{( \zeta^2 - k^2)^2}$]. 
 There is a constant $C$,  such that for $ | \uE | \le 1 $,  $k \ge \max\{ 1,  4 | \uE |^2 \} $, 
  $t \in [0, 1]$ and $ 2 | \uE |^2  z \le 1$ there holds 
  \begin{equation}
\label{2.28}
\big| \mathcal{E}_1 (t, z ) \big|  \le C \frac{ k    \ e^\rho}{ \sqrt {1+ \rho}} 
\end{equation}
\begin{equation}
\label{2.29}
\big| \mathcal{E}_2 (t, z ) \big|  \le C \frac{ k^2     \ e^\rho}{  \tilde  r  \sqrt {1+ \rho}} . 
\end{equation}
with $\frac{1}{\tilde r}  =   \frac{1}{k} + \frac{ t}{1 + \rho} $. 
\end{prop}

\begin{proof}  For $\cE_1$  we proceed as for    $\cE_0$,  noticing that 
$| \frac{k^2}{ \zeta^2 - k^2} | \approx \frac{ k}{ r} $ on the integration path. 

For \eqref{2.29}, we make the same choice of radius $  r$ when 
$k | E|^2 z t \ge 2$. Noticing that 
$| \frac{k^4}{ (\zeta^2 - k^2)^2 } | \approx \frac{ k^2}{   r^2} $ on the integration path 
yields the estimate
\begin{equation*}
\big| \mathcal{E}_2 (t, z ) \big|  \le C   \sqrt { \frac{t}{ | \uE |^2 kz}} \
 \frac{k^2  e^\rho } {\sqrt{1+   \rho}  } . 
\end{equation*}

When $k | E|^2 z t \le 2$, 
the   length $\Gamma $  is $O(r)$, we maximize $r$ in the range $[1, \mez k]$ 
instead of minimizing it as we did in  
 \eqref{2.26},  so that the imaginary parts of the phases remain bounded. Namely, we choose 
\begin{equation*}
\tilde r = \min \{ \frac{1}{t} , \mez k \} 
\end{equation*}
  and get 
\begin{equation*}
\big| \mathcal{E}_2 (t, z ) \big|  \le C \frac{k^2}{\tilde  r } . 
\end{equation*}
in this zone. Combining the estimates yields \eqref{2.29}
\end{proof}

\medbreak

\subsection{Integral formulas for the solutions of \eqref{2.4} \eqref{2.5}} 

At the given frequency $k$, the system \eqref{2.4} \eqref{2.5} is a classical hyperbolic 
Goursat problem in dimension 2, and therefore has a unique smooth solution for smooth data
(see \cite{Al}), which respects the causality principle in time. 
We will use the following description: 

\begin{prop}
\label{prop2.7}
For $F = {}^t (f, g)$ and $h$ in $C^\infty ([0, T] \times [0, Z])$ and 
the problem \eqref{2.4} \eqref{2.5}  with $U_0 = {}^t (u_0, v_0) = 0 $, 
has a unique solution $U = {}^t(u, v)$, $n$ in $C^\infty ([0, T] \times [0, Z])$ and 
\begin{equation}
\label{2.27a}
\begin{aligned}
U (t, z) = &   i   \int_0^z e^{ i (z - z') A(\infty) }  F( t, z' ) dz'        
\\
& +  \int_0^z \!\!\!\! \int_0^t \cE_0(t',  z') F(t- t', z- z') dz' dt' 
\\
 & +  \int_0^z \!\!\!\! \int_0^t \cE_1(t',  z') \begin{pmatrix}  -i\\ i \end{pmatrix}
 h  (t- t', z- z') dz' dt'  , 
\end{aligned}
\end{equation}
\begin{equation}
\label{2.27b}
\begin{aligned}
n (t, z) = &   
   \int_0^z \!\!\!\! \int_0^t  \begin{pmatrix}  -i, & i \end{pmatrix} \cE_1(t',  z') 
 F  (t- t', z- z') dz' dt'       
\\
& +  i \int_0^z \!\!\!\! \int_0^t \cE_2 (t',  z') h(t- t', z- z') dz' dt' . 
\end{aligned}
\end{equation}
\end{prop}


\subsection{Exponentially growing solutions} 

We show that the rate of amplification $e^\rho$ observed in 
Propositions~\ref{prop2.5} and \ref{prop2.6} is sharp. 
Consider the solution of 
 \begin{equation}
\label{2.30}
\left\{ \begin{aligned}
& i \D_z  u  + \Delta_x   u -    \underline E   n =  0  ,\\
& i \D_z  v  -  \Delta_x   v  +     \underline E   n =    0   , \\ 
& (\D_t^2 - \Delta_x)  n -  \underline E   \Delta_x  ( u + v  ) =  0. 
\end{aligned}
\right.
\end{equation} with initial-boundary conditions
\begin{equation}
\label{2.31} 
n_{| t = 0} = \D_t n _{| t = 0} = 0, \quad  u_{| z = 0} =  
\frac{\sin (k t)}{k}  e^{ i \xi  x} ,  \quad  v_{| z = 0} = 0. 
\end{equation}
 This amounts  to solve  
  \begin{equation}
\label{2.32}
\left\{ \begin{aligned}
& i \D_z  u  - k^2  u -    \underline E   n = 0  ,\\
& i \D_z  v  + k^2    v  +     \underline E   n =    0   , \\ 
& (\D_t^2 - \Delta_x)  n -  \underline E   k^2   ( u + v  ) = 0  . 
\end{aligned}
\right.
\end{equation} 
with $k = | \xi |$, together with initial-boundary conditions  
\begin{equation}
\label{2.33} 
n_{| t = 0} = \D_t n _{| t = 0} = 0, \quad  u_{| z = 0} =  \frac{\sin (k t)}{k}    ,  \quad  v_{| z = 0} = 0. 
\end{equation}
Because the Fourier-Laplace transform of   $1_{ \{ t > 0 \}} \frac{\sin (kt)}{k} $ 
is $\frac{1}{ k^2 -  \zeta^2 } $, the solution is 

  \begin{equation}
 \label{2.34} 
U (t, z) = \begin{pmatrix} u(t, z) \\ v(t, z)  \end{pmatrix}
 =   \frac{1}{2 i  \pi} 
\int_{\RR - i \gamma}  e^{  i \big( t \zeta + z A_k (\zeta)  \big)  } R     \frac{ d \zeta}{k^2 -\zeta^2} , 
\end{equation}
 \begin{equation}
 \label{2.35} 
n  (t, z) 
 =  -   \frac{1}{2 i  \pi} 
\int_{\RR - i \gamma}  L  \  e^{  i \big( t \zeta + z A_k (\zeta)  \big)  } R  \ 
 \frac{ k^2 E } {( \zeta^2  -  k^2)^2} \   \frac{ d \zeta}{\zeta}  ,  
\end{equation}
where 
$$
R = \begin{pmatrix}
    1     \\
  0    
\end{pmatrix},\qquad 
L = \begin{pmatrix}
   1   \ , \  1
\end{pmatrix}. 
$$
      
 \begin{theo}
 \label{th2.7} 
There are constants  $C \ge c > 0$ such that for all   $k \ge 2 $, $t \in [0, 1]$, 
$z \in [0, 1] $  :  
 \begin{equation*}
 \begin{aligned} 
  \big|   U(t, z )\big |  & \le   C \frac{1}{ k}  
   \frac{ e^{    \rho  } }{\sqrt {1 + \rho} }    
\\ 
\big| n (z, t) \big|  &    \le C  \big(  \frac{1}{ k}  + \frac{t}{1 + \rho}   \big)  \frac{ e^{    \rho  } }{\sqrt {1 + \rho} }  .  
\end{aligned}
\end{equation*}  
Moreover, for $\rho \ge 1$, $t \ge \mez$   and $k$ large enough,  
\begin{equation}
\label{2.36}
  \big| n_k(t, z) \big|  \ge c    \rho ^{- \frac{5}{2}}   e^{    \rho  } \Big( 1 + O( \rho^{-1} )\Big) . 
  \end{equation}

 \end{theo} 
 
   \begin{proof}   One can deform the integration path $\RR - i \gamma$ 
   to   a contour $\Gamma$ which is   the union of the circles 
$$
\Gamma_1 = \{  | \zeta - k | = r \} , \qquad  
\Gamma_2  = \{  | \zeta + k | = r \}  
$$
with $ r  \in [1, \mez k] $.   
The upper bounds follow from Proposition~\ref{prop2.6}.

   We now concentrate on the proof of the lower bound for $n$
assuming that $\rho = \sqrt{ 2 | \underline E |^2 z t}  \ge  2 $.  
We choose the radius $r = \frac{\rho}{2t} $ as in \eqref{2.22}.  
By \eqref{2.17}, for $j = 1, 2$,  the integral over $\Gamma_j $  contributes to  
$$ 
\begin{aligned}
N_j   (t, z) 
 =   &  \frac{1}{2 i  \pi} 
\int_{\Gamma_j }    \  e^{  i   t \zeta + i z k \lambda}   \frac{ \lambda -  k} {2  \lambda } \ 
 \frac{ k^2 E } {( \zeta^2  -  k^2)^2 } \   \frac{ d \zeta}{\zeta}
 \\
&  +   \frac{1}{2 i  \pi} 
\int_{\Gamma_j }    \  e^{  i   t \zeta -   i z k \lambda}   \frac{ \lambda +   k} {2  \lambda } \ 
 \frac{ k^2 E } {( \zeta^2  -  k^2)^2} \   \frac{ d \zeta}{\zeta} . 
\end{aligned}
 $$ 
 Consider first $N_1$. 
  We evaluate the phases $\vp = t \zeta  - z k  \lambda$  and 
 $\psi = t \zeta  +  z k  \lambda$ on  $\Gamma_1$ as in \eqref{2.23} \eqref{2.24}.  For 
$\zeta = k + r  e^{ i \theta}$ we have 
$$
\begin{aligned}
\vp ( k + r  e^{ i \theta})   & =  k t  - k^2 z +   2 i 
\rho  \sin \theta  -  \sigma_k( r, \theta), 
\\
\psi  ( k + r  e^{ i \theta})   & =  k t  +  k^2 z +  2 
\rho \cos \theta + \sigma_k (r,  \theta). 
\end{aligned}
$$
where  the $\sigma_k$ are smooth functions of  $\theta \in \TT$, uniformly 
bounded as well as their derivatives when  $1 \le r \le \mez k$.

 Note that the imaginary part of $\psi$ is bounded,  while by \eqref{2.20} 
\begin{equation}
\label{2.37}
\Big|  \frac{ \lambda -  k} {2  \lambda } \ 
 \frac{ k^2 |\uE|^2 } {  ( \zeta^2  -  k^2)^2 } \Big|  \le  \frac{C | \uE |^2}{ k  r^3} \le \frac { C | \uE |^2 }{k} 
\end{equation}
on $\Gamma_1$. Thus the corresponding term contributes to $O(k^{-1})$  terms 
in $N_1$. 
It remains to evaluate 
$$
\tilde N_1 (z, t)  = 
   \frac{1}{2 i  \pi} 
\int_{\Gamma_1 }    \  e^{  i   t \zeta -   i z k \lambda}   \frac{ \lambda +   k} {2  \lambda } \ 
 \frac{ k^2 E } { (\zeta^2  -  k^2)^2} \  { d \zeta} 
$$
We have the following expansion 
$$
 \frac{ \lambda +   k} {2  \lambda } \ 
 \frac{ k^2 E } {  ( \zeta^2  -  k^2)^2 }  =  \frac{ e^{ - 2  i \theta}} {r^2 } \  
  p_k (r, \theta)  
$$
where  the $p_k$ are smooth functions of  $\theta \in \TT$, uniformly 
bounded as well as their derivatives when  $1 \le r \le \mez k$ and such that, uniformly
on this domain
$$
\lim_{ k \to \infty}  p_k(r, \theta )  =  \frac{  \uE } {2} 
$$
Substituting,  we see that
 $$
 \tilde N_1 (z, t)   =  \frac{ e^{ i (k t - k^2 z) }}{2 \pi r^2}  \int_0^{2 \pi} 
   e^{ -  2 \rho  \sin \theta  }    \   p_k (r,  \theta)   e^{-  i z \sigma_k (r, \theta) - 2 i \theta}    d \theta  . 
$$
The  stationary phase theorem to  implies that 
 $$
 \tilde N_1 (z, t)   =   e^{ i (k t - k^2 z) } \ \frac{ e^{ 2 \rho }} {2 r^2 \sqrt {\pi \rho } } 
  \big( p_k (r, - \pi/2) e^{ - i \sigma_k( r, - \pi/2) + i \pi }   +  O( \rho^{-1}) \big) 
 $$ 
 Therefore, for $k$ large enough, 
 \begin{equation*}
  \big| N_1 (t, z) \big|  \ge  \frac{ | \underline E | t^2 }{  \sqrt \pi   \rho^{5/2}  }    e^{    \rho  } 
   \Big( 1 + O( \rho^{-1} )\Big)  + O(1) . 
  \end{equation*}
  
  \medbreak
  {\bf b) } The analysis of the contribution of the integral over $\Gamma_2$ is similar,  
  except that the  phase with large imaginary part  is now $\psi$. 
  Noticing that the corresponding amplitude is bounded by \eqref{2.37}, we see that 
  \begin{equation*}
  \big| N_2 (t, z) \big| \le \frac{C }{r^2}    \frac{ e^{ \rho}} {\sqrt \rho} + O(1) 
  \end{equation*} 
  and because $\frac{1}{r} \approx  \frac{t}{\rho}  \lesssim \frac{1}{\rho}$, 
   the estimate \eqref{2.36} follows for $k$ and $\rho$ large enough.
\end{proof}


\section{The nonlinear instability}

\subsection{The method of proof}

We fix a constant field $\uE \ne 0$, noticing that 
$(\uE, 0)$ is a solution of  \eqref{1.3}.  We compare it to another solution  of the form 
\begin{equation} 
\label{3.1}
E  = \underline E  + e^a  + \tilde e, \quad   n = n^a + \tilde n 
\end{equation} 
where   $(e^a,  n^a)$ is a  solution of the homogeneous linearized 
equation \eqref{2.1} \eqref{2.2},  with small boundary value on $\{ z = 0 \}$ and 
such that $n^a$ and $n$ are arbitrarily large at an arbitrarily small distance $z > 0$. 
Since we consider periodic functions in $x$, 
it is sufficient  to exhibit an example where  \emph{ $x$ is one dimensional}, which we now assume. 
We also assume, without restriction,   that \emph{ $\uE$ is real.}

We choose $(e^a, n^a)$ to be the solution of 
\begin{equation}
\label{a3.1}
\left\{ \begin{aligned}
& i \D_z  e^a  + \Delta_x   e^a -    \uE   n^a = 0  ,\\
& (\D_t^2 - \Delta_x)  n^a -  \Delta_x (    \uE   e^a + \uE \overline {e^a} ) =  0  
\end{aligned}
\right.
\end{equation}
\begin{equation}
\label{3a.2}
n^a _{| t = 0}  =  \D_t n^a_{| t = 0} = 0, \quad  e^a_{ | z = 0}  = \alpha   \frac{\sin (kt)}{k}  e^{ i k x} ,  
\end{equation}
where the frequency $k$ is large and $\alpha  $ is a small parameter. Typically, we will choose 
$\alpha  \approx  k^{ - \sigma }$ so that the boundary data  is of order
$k^{ s- \sigma}$ in $H^s$. Thus
 \begin{eqnarray}
 &&e^a (t,z, x) = u(t, z) e^{ i k x}  + \overline{ v (t, z)} e^{ - i k x}  , \\
 &&n^a (t,z, x) = \nu (t, z) e^{ i k x}  + \overline { \nu (t, z)} e^{ - i k x} 
 \end{eqnarray} 
where $(u^a, v^a, \nu^a )$ solve 
 \begin{equation}
\label{a3.6}
\left\{ \begin{aligned}
& i \D_z  u^a    - k^2    u^a   -    \uE   \nu^a =  0   ,\\
& i \D_z  v^a   +  k^2    v^a     +      \uE     \nu^a =   0   , \\ 
& (\D_t^2 + k^2 )  \nu^a + k^2    \uE   ( u^a  + v^a  ) = 0, 
\end{aligned}
\right.
\end{equation} 
 \begin{equation}
\label{a3.7}
\nu^a _{| t = 0}  =  \D_t  \nu^a_{| t = 0} = 0, \quad  u^a_{ | z = 0}  = \alpha \frac{\sin(kt)}{k} , \quad 
 v^a_{ | z = 0}  = 0 .
\end{equation}

  The solution $(e^a, n^a)$ is amplified as described in Theorem~\ref{th2.7}. 
  We construct correctors  $(\tilde e, \tilde n)$ such that 
  $(E, n)$ given by \eqref{3.1} is an exact solution of the nonlinear system, 
  on a small domain $\Omega$ where  $| \tilde n|  \ll  |n^a | $ and such that $n^a$ and thus $n$
  is arbitrarily large on a part of the boundary $\D \Omega$. 
The equations for the remainders $(\tilde e, \tilde n)$ in \eqref{3.1} read
\begin{equation}
\label{6.2}
\left\{ \begin{aligned}
& i \D_z  \tilde e  + \Delta_x   \tilde e -    \underline E   \tilde n =    
(e^a + \tilde e) (n^a + \tilde n)  ,\\
& (\D_t^2 - \Delta_x)  \tilde n -  \Delta_x (   \overline {\underline E}  \tilde e + \underline  E \overline {\tilde e } ) =   \Delta_x \big( ( e^a + \tilde e) (\overline {e^a} + \overline {\tilde e}) \big). 
\end{aligned}
\right.
\end{equation}
We add homogeneous  initial-boundary data  
 \begin{equation}
\label{a3.9}
\tilde n _{| t = 0}  =   \D_t \tilde n _{| t = 0} = 0, \quad  \tilde e_{ | z = 0}  = 0 . 
\end{equation}
so that  $E$ has the same boundary data as $ \uE + e^a$.

We solve the nonlinear equation \eqref{6.2} by Picard's iterations, in suitable Banach spaces
of analytic functions which we now describe.


\subsection{Construction of the correctors  } 

The approximate solution $(e^a, n^a)$ has only two frequencies in $x$, $+k$ and $-k$. 
The nonlinear interaction will create all the harmonics. This leads to consider 
functions of the form
\begin{equation} 
\label{a5.1}
 u   = \sum_{p \in \ZZ}  u_p (t, z)  e^{ i k p  x} . 
\end{equation}

\begin{defi}

Given parameters $k \ge 2 $, $\delta \in ]0, 1]$, $b \in ]0, 1[$ 
  and 
$s \ge 1$, we denote by  $\EE $  the space of  functions 
\eqref{a5.1} 
such that 
\begin{equation}
\label{a5.2}
| u_p(t, z) | \le  C    \  \frac{1}{ (1 + | p |)^s }  \ \delta^{\la p \ra } 
  e^{ \la p \ra \rho} 
\end{equation}
with 
\begin{equation}
\la p \ra = \max \{ 2 , | p | \}, \qquad \rho = \sqrt{ 2 k | \uE |^2 z t } , 
\end{equation}
on the domain 
\begin{equation}
\Omega   = \big\{ (t, z, x) \in [0, 1] \times [0, b] \times \TT \ :  
\   \delta e^{\rho}  \le  1  \big\}. 
\end{equation}
The best constant in \eqref{a5.2} is the norm of $u$ in $\EE$. 
\end{defi} 

\begin{rem}
\textup{The factor  $\delta e^\rho \le 1 $ present for the frequencies $\pm k$ in the approximate 
solution is expected to intervene at the power $| p |$    in the harmonic $ pk$ by successive
nonlinear interaction.  Note that  the series \eqref{a5.1} define real analytic functions in $x$ 
when    $(t, z) \in \{ \delta e^\rho  < 1\} $.   }
\end{rem}

The space $\EE$ and the domain $\Omega$ depend  strongly on the parameters
$k, b$ and $\delta$. However, to lighten notations, we do not mention this  dependence explicitly.

\begin{lem}
\label{lema3.2}
 Consider the solution $(e^a, n^a) $ of \eqref{a3.1} \eqref{3a.2}  with 
$\alpha = \gamma \delta$.  

 There is a constant  $C_0$  independent of $k \ge 2$, $b\in ]0, 1] $, 
 $\gamma > 0$  and $\delta \in ]0, 1]$ such that  

\quad i)  
$f^a = n^a u^a  \in \EE $  and $h^a = u^a   u^a  \in \EE $ and 
\begin{equation}
\|  f^a \|_{\EE}  \le   C_0     \frac{  \gamma^2 }{k}     , \quad 
\|  h^a \|_{\EE }  \le   C_0      \frac{\gamma^2 }{k^2}  . 
\end{equation}

\quad ii) for all  $v \in \EE$,   $u^a v$   and $n^a v  $   belong to $\EE$ and 
\begin{equation}
\|  u^a v   \|_{\EE}  \le  C_0     \frac{  \gamma   }{k}    \| v \|_{\EE}   , \quad 
\|  n^a  v  \|_{\EE }  \le   C_0  \gamma \| v \|_{\EE}. 
\end{equation}
\end{lem}

\begin{lem}
\label{lema3.3} $\EE$ is a Banach algebra and there
 is a constant  $C_0$  independent of $k$, $b$  and $\delta$ such that for all 
 $u $ and $v$ in $\EE$:
 \begin{equation}
 \| u v \|_{\EE}   \le   C_0 \| u \|_{\EE} \ \| v \|_{\EE} . 
 \end{equation}
\end{lem}

\begin{cor}
\label{corn3.5}
 Introduce the notations
\begin{eqnarray}
&& \Phi (e, n) =  (e^a + e) (n^a + n) 
\\
&&
\Psi (u, v) :=  (e^a + u) (\overline {e^a} + \overline v) . 
\end{eqnarray}
Then 
\begin{equation*}
\|  \Phi (e, n)  \|_{\EE}  \le   C_0   \Big(   \frac{  \gamma^2 }{k}      + 
 \frac{\gamma}{k}   \| n \|_{\EE}
 +  \gamma \| e  \|_{\EE}  \big) +  \| e \|_{\EE} \ \| n \|_{\EE}\Big). 
\end{equation*}
\begin{equation*}
\|  \Psi (u, v)   \|_{\EE}  \le   C_0   \Big(     \frac{  \gamma^2 }{k^2}   + 
 \frac{\gamma}{k} \big(  \| u  \|_{\EE}
 +   \| v \|_{\EE}  \big) +  \|  u  \|_{\EE} \ \| v \|_{\EE}\Big). 
\end{equation*}

\end{cor}

\bigbreak

Next, consider the linear problem
\begin{equation}
\label{a3.17}
\left\{ \begin{aligned}
& i \D_z    e  + \Delta_x     e -    \underline E     n =   f   ,\\
& (\D_t^2 - \Delta_x)   n -  \Delta_x (   \overline {\underline E}    e + \underline  E \overline {  e } ) =  \Delta_x h  , 
\\
&  n _{| t = 0}  =   \D_t   n _{| t = 0} = 0, \quad   e_{ | z = 0}  = 0 . 
\end{aligned}
\right.
\end{equation}
 
\begin{prop}
\label{prop35} 
For $f   $ and $h $ in $\EE$, the solution $(e, n)$  of \eqref{a3.17}    belongs to $\EE$ and 
\begin{eqnarray}
&&
\| e \|_{\EE}  \le  C_1 \,  b  \,  \| f \|_{\EE}  +  C_1  \ln k  \|  h \|_{\EE}
\\
&&
\| n \|_{\EE}  \le  C_1     \frac{1}{| \uE |^2}  \ln  ( | \ln \delta | )   \| f \|_{\EE}  +   C_1 
\frac{ k}{| \uE|^2}   \|  h \|_{\EE}
\end{eqnarray}
where $C_1$  is  independent of $k$, $b$  and $\delta$
\end{prop}

Denote by $ (e ,n)  = \cT(f, h)$  the solution of \eqref{a3.17}. The equation \eqref{6.2} reads
\begin{equation}
\label{6.3}
(\tilde e, \tilde n)  = \cT \big( \Phi (\tilde e, \tilde n) ,  \Psi (\tilde e, \tilde e )  \big) 
:=  \cF (\tilde e, \tilde n) := \big( \cE(\tilde e, \tilde n) , \cN(\tilde e, \tilde n) \big) . 
\end{equation}

 \begin{prop}
 For $\gamma \ge 1$ and 
 \begin{eqnarray}
 \label{n3.23}
&&   4 C_0 C_1 \gamma^2  \big(  b  +  \frac{\ln k}{k}\big)  \le 1, 
 \\
\label{n3.24}
 &&  4  C_1 C_0 \gamma^2   
   \big(    \frac{\ln ( \ln \delta |) }{k}  +  \frac{1}{k}  \big) \le | \uE |^2 , 
 \end{eqnarray}
 $\cT$ maps  $\BB_{\frac{1}{k}} \times \BB_1$ into itself, 
 where $\BB_R$ denotes the ball of radius $R$ in $\EE$, and is a contraction for the 
 norm
 $$
 \| e \|_\EE + \frac{1}{k} \| n \|_{\EE}. 
 $$ 
 if $\gamma > 2$. 
 \end{prop}
 
 \begin{proof} 
 By Corollary~\ref{corn3.5},  for $(e, n) \in \BB_{\frac{1}{k}} \times \BB_1$ and 
 $\gamma  \ge 1$, there holds
 $$
 \| \Phi (e, n) \|_{\EE}  \le  4 C_0 \frac{\gamma^2}{k}, 
 \quad   \| \Psi (e, e) \|_{\EE}  \le  4 C_0 \frac{\gamma^2}{k^2}.  
 $$
 Thus, by Proposition~\ref{prop35} and \eqref{n3.23} \eqref{n3.24}
 $$
 \| \cE (e, n) \|_{\EE}  \le  4  C_1 C_0 \big(  b  \frac{\gamma^2}{k}  + 
  \frac{\gamma^2 \ln k }{k^2}\big)  \le \frac{1}{k}  
 $$
 and 
  $$
 \| \cN  (e, n) \|_{\EE}  \le  \frac{4  C_1 C_0 \gamma^2 }{|\uE|^2} 
  \big(    \frac{\ln ( \ln \delta |) }{k}  +  \frac{1}{k}  \big) \le 1. 
 $$
 This shows that  $\cT$ maps  $\BB_{\frac{1}{k}} \times \BB_1$ into itself.

Consider $(e, n) $ and $(e', n') $ in $\BB_{\frac{1}{k}} \times \BB_1$. Denote by 
$\delta e = e- e'$, $\delta n = n- n'$, $\delta \Phi = \Phi(e, n) - \Phi(e', n')$ etc. 
There holds
\begin{eqnarray*}
&&\|  \delta \Phi \|_{\EE}  \le  2 C_0 \gamma  \big(  \| \delta e \|_\EE + \frac{1}{k} \| \delta n \|_{\EE} \big), 
\\
&&\|  \delta \Psi \|_{\EE}  \le \frac{ 4  C_0 \gamma}{k}    \| \delta e \|_\EE , 
\end{eqnarray*}
and thus, using Proposition~\ref{prop35} and \eqref{n3.23} \eqref{n3.24}, 
\begin{eqnarray*}
&&\|  \delta \cE  \|_{\EE}  \le  \frac{1}{\gamma}   \big(  \| \delta e \|_\EE + \frac{1}{k} \| \delta n \|_{\EE} \big), 
\\
&&\frac{1}{k}  \|  \delta \cN  \|_{\EE}  \le \frac{1}{ \gamma}  \big(  \| \delta e \|_\EE + \frac{1}{k} \| \delta n \|_{\EE}  
\end{eqnarray*}
implying that $\cT$ is a contraction if $\gamma > 2$. 
 \end{proof}

\begin{cor}
  Assume that  $\gamma > 2$ and that  \eqref{n3.23} and \eqref{n3.24} are satisfied. 
  Then   the equation \eqref{6.2} has a unique solution in 
   $\BB_{\frac{1}{k}} \times \BB_1$.

\end{cor}

\subsection{Choice of parameters, proof of Theorem~\ref{mainth}} 

Fix  $s = 1$,  $\epsilon \in ]0, 1[ $,  $\epsilon' \in ]0,  \mez  \epsilon [ $ and $\sigma  > 0$. We choose
\begin{equation}
\label{6.10}
\delta = k ^{ - \sigma}, \quad  b = k^{ - \epsilon} , \quad  \gamma = k^{\epsilon'}.
\end{equation}
Then the conditions \eqref{n3.23} and \eqref{n3.24} are satisfied for $k$ large
and there is a solution  
 \begin{equation*} 
E  = \underline E  + e^a  + \tilde e, \quad   n = n^a + \tilde n 
\end{equation*} 
of the original system, with 
\begin{equation}
\| \tilde e \|_{\EE} \le  \frac{1}{k}, \qquad \| \tilde n \|_{\EE} \le 1. 
\end{equation}
 Therefore the first Fourier coefficient of $n$ satisfies 
\begin{equation}
\label{n3.27}
| n_1 (t, z) - n^a_1 (t, z)  |  \le   \delta^2 e^{ 2 \rho}   \le 1 
\end{equation}
while, by Theorem~\ref{th2.7}, for $t \ge \mez$ and $\rho $ large enough 
\begin{equation}
| n^a_1 (t, z) | \ge c \gamma  \rho^{-5/2}  \delta e^{ \rho} . 
\end{equation}
 
Note that $\delta e^{ \rho} = 1 $ for  $\rho = \sqrt { 2 k | E |^2 zt } = \sigma \ln k $ 
thus for $z t = \alpha k^{-1} (\ln k)^2$ for some constant $\alpha > 0$. 
For $k$ large,  $ \alpha k^{-1} (\ln k)^2 \ll k^{- \epsilon}$ and therefore  
\begin{equation}
 \Gamma := \overline \Omega \cap \{  \delta e^{ \rho} = 1  \}  \ne \emptyset . 
\end{equation}
On this set,   
\begin{equation}
| n^a_1 (t, z) | \ge c'    \frac{  k^{ \epsilon'} }{   (\ln k)^\frac{5}{2}} . 
\end{equation}
 which tends to $+ \infty$ with $k$. 
 With \eqref{n3.27}, this finishes the proof of Theorem~\ref{mainth}.


\section{Proofs}

\subsection{Lemmas}

We collect here several elementary estimates which will be 
used repeatedly in the sequel.

\begin{lem}
\label{lem3.1} There is a constant $C$  such that for all 
$\lambda > 0$,   $t \ge 0$ and all $z \ge 0$: 
\begin{equation}
\label{aa3.1}
\int_0^z \int_0^t  e^{ - \sqrt{ \lambda z' t'}} dz' dt'  \le   C  
\frac{1}{\lambda}  \ \ln \big( 1 + \lambda z t \big),  
\end{equation}
\begin{equation}
\label{aa3.2}
\int_0^z \int_0^t   z'    e^{ - \sqrt{ \lambda z' t'}} dz' dt'  \le 
\frac{ C  z }{\lambda}    , 
\end{equation}
\begin{equation}
\label{aa3.3}
\int_0^z \int_0^t   t'    e^{ - \sqrt{ \lambda z' t'}} dz' dt'  \le  
\frac{ C  t }{\lambda}   . 
\end{equation}

\end{lem}

\begin{proof}
$$
\int_0^t  e^{ - \sqrt{ \lambda z' t'}}  dt'  =  \frac{1}{\lambda z'}   H( \lambda z' t) 
$$
with 
$$
H(u) = \int_0^u e^{ - \sqrt {u'}} du' \le  C    \frac{u}{ 1+ u} . 
$$
Thus 
$$
\int_0^z \int_0^t  e^{ - \sqrt{ \lambda z' t'}} dz' dt'  \le C  
\int_0^z \frac{t}{1 + \lambda z' t} \ dz' 
= \frac{C }{\lambda} \int_0^{\lambda z t}   \frac {ds }{1 + s}. 
$$
and \eqref{aa3.1} follows. The other two estimates are symmetric in $t$ and $z$. 
Integrating in $t'$ first as above we see that  
$$
\int_0^z \int_0^t  z'   e^{ - \sqrt{ \lambda z' t'}} dz' dt'  \le   C 
\int_0^z \frac{ t z' }{1 + \lambda z' t} \ dz' 
\le  \frac{C }{\lambda} \int_0^{z}   dz'   
$$
implying \eqref{aa3.2}. 
\end{proof}

\begin{lem}
\label{lem3.1b}
For $p \ge 1$ and $\mu \ge 8$, 
\begin{equation*}
\frac{1}{p} \ln ( 1 + \mu p^2 ) \le   \ln (1 + \mu). 
\end{equation*}
\end{lem}
\begin{proof}
The upper bound of the left hand side for real $p \ge 1$  is equal to 
$ c_0  \sqrt \mu $  if   $\mu \le y_0^2$ with  $ y_0$ being  the positive root of 
$\ln (1 + y_0^2 ) =  \frac{ 2 y_0^2}{1 + y_0^2} $, where the function 
$\frac{\ln (1+ y^2)}{ y}$ reached its maximum equal to $c_0 = \frac{2 y_0}{ 1 + y_0^2}$. 
When  $\mu \ge y_0^2$,  the upper bound is attained at $p =1$ and thus equal 
to $\ln (1 + \mu)$. Because $y_0^2 \le 8$, the lemma follows. 
\end{proof}


\subsection{Linear estimates} 
In this section we prove the estimates of Proposition~\ref{prop35}.  
Expanding  the system \eqref{a3.17} in Fourier series in $x$ and 
denoting by $(e_p, n_p, f_p)$ and $h_p$ the Fourier coefficients, yields the following system 
for $U_p :=  (u_p, v_p) =   (e_p, \overline{e_{-p}} )$
 \begin{equation}
\label{b4.1}
\left\{ \begin{aligned}
& i \D_z  u_p   - k^2 p^2    u_p   -    \underline E   n_p  =  f_p    ,\\
& i \D_z  v_p   +  k^2 p^2     v_p   +     \underline E    n_p =   g_p   , \\ 
& (\D_t^2 + k^2 )  n_p + k^2 p^2   \underline  E   ( u_p  + v_p  ) = - k^2p^2  h_p , 
\end{aligned}
\right.
\end{equation}
with $g_p = - \overline{f_{-p}}$,  plus initial boundary boundary conditions:   
 \begin{equation}
\label{b4.2}
n _p{}_{| t = 0}  = 0, \quad \D_t n_p{}_{| t = 0} = 0, \quad  u_p{}_{ | z = 0}  =  
 v_p{}_{ | z = 0}  = 0 .
\end{equation}
Proposition~\ref{prop35} is an immediate corollary of the following estimates, where we use the notations 
\begin{equation}
\label{bb4.2}
\rho(t, z) = \sqrt{ 2 k | E |^2 z t } . 
\end{equation}
\begin{prop}
\label{propb4.1} 
There is a constant $C_1$  is  independent of $k  \ge 2 $, $b\le 1 $  and $p $ such that for 
$F_p := (f_p, g_p)$ and $h_p$ satisfying on $\Omega$ 
\begin{equation}
\label{b4.3}
| F_p(t, z) | \le  A   e^{ \la p \ra  \rho(t, z) } , \quad   | h_p(t, z) | \le  B   e^{ \la p \ra  \rho(t, z) } ,
\end{equation}
 the solution $(U_p, n_p)$  of \eqref{b4.1} \eqref{b4.2} on $\Omega$
satisfies  
\begin{eqnarray}
\label{b4.8}&&
| U_p(t, z) |   \le  C_1  \Big( b  \,  A  +      \frac{1}{| \uE |^2}    \ln ( \ln \delta |  )\,  B   \Big)
 e^{ \la p \ra  \rho(t, z)}, 
\\
\label{b4.9}&&
|  n_p  (t, z) |  \le  C_1 \Big(  \frac{1}{ | \uE |^2 }  \ln ( |  \ln \delta |) A     +     
 \frac{1}{ | \uE |^2 } \big( k  +    \ln ( |  \ln \delta |) \big)  B \Big) e^{ \la p \ra  \rho(t, z)}  . 
\end{eqnarray}
 
\end{prop}

\begin{proof}    We use Proposition~\ref{prop2.7} at the frequency $kp$. 

{\bf a) }   By \eqref{2.27a},  $U_p $ is the sum of three  terms.  The first one is bounded by 
\begin{equation}
\label{b4.10}
A \int_0^z  e^{ \la p \ra  \rho(t,z' ) }  dz'   \le A z e^{ \la p \ra  \rho(t, z)}   . 
\end{equation}
When $p= 0$, only this term is present. When $p \ne 0$, 
 the second term is  
\begin{equation*}
 U_{p, 0}   :=  \cE_0 * F_p   = \int_0^z \!\!\!\! \int_0^t \cE_0(t',  z') F_p (t- t', z- z') dz' dt' .
\end{equation*}
Note the following identity for non negative real numbers:  
\begin{equation*}
 \sqrt{ a' b'}  + \sqrt {a'' b''}  \le \sqrt{ (a'+a'') (b'+b'')} . 
\end{equation*}
In particular, 
\begin{equation}
\label{a4.6}
\rho (t', z')  + \rho (t-t', z-z')  \le \rho (t, z), 
\end{equation}
and therefore, by Proposition~\ref{prop2.5}  and \eqref{b4.3}: 
\begin{equation*}
| U_{p,0} (t, z )|  \lesssim \alpha   e^{\la p \ra \rho(t, z)  }  
 \int_0^z \!\!\!\! \int_0^t   e^{ - \la  p \ra  \rho(t', z') }    
 \frac{ r'_p    e^{ \rho'_p } }{\sqrt{(1 + \rho'_p)} }   dz' dt'  
\end{equation*} 
where  $\rho'_p  = \sqrt{ 2 k p  |E|^2 z' t' }$ and 
$r'_p = 1 + \frac { kp | \uE|^2  z'}{1 + \rho'_p}\le   1+  kp | \uE|^2 z' $ . 

When $p \ge 2$, 
\begin{equation*} 
 p \rho(t', z')  - \rho'_p = \sqrt{ 2 (p^2 - p) k |E|^2 z' t'}  \ge      \sqrt{   k p^2  |E|^2 z' t' }
 \ge \mez p \rho(t', z'). 
\end{equation*} 
When $p = 1$, 
\begin{equation*}
 \la 1 \ra  \rho(t', z')  - \rho'_1  =   2  \rho(t', z') -   \rho(t', z') =   \rho(t', z'). 
\end{equation*} 
Thus, in any case, 
\begin{equation}
\label{a4.8} 
\la  p \ra  \rho(t', z')  - \rho'_p    \ge       \mez p \rho(t', z'). 
\end{equation} 

Therefore
\begin{equation*}
| U_{p,0}  (t, z )|  \lesssim A    e^{\la p \ra  \rho (t, z)  }  
 \int_0^z \!\!\!\! \int_0^t   (1 + k p | \uE |^2  z')   e^{ - \mez  p \rho(t', z') }    
    dz' dt' . 
\end{equation*} 
Using Lemmas~\ref{lem3.1}  with $\lambda  = \mez  k p^2 | \uE |^2 $,  yields 
\begin{equation}
\label{b4.13}
| U_{p, 0}  (t, z )|  \lesssim A  \Big(  \frac{1}{ \lambda  }   \ln \big(1 + \lambda   z t \big)   +  
kp | \uE |^2 \frac{ z }{ \lambda}   
\Big)  e^{p \rho }    \lesssim A z e^{ \la p \ra \rho },  
\end{equation}
where we have used that $ \frac{1}{ \lambda  }   \ln \big(1 + \lambda   z t \big)  \le zt $. 
\medbreak

Similarly  $U_{p, 1}  = \cE_1 * h_p $ satisfies 
\begin{equation*}
\begin{aligned}
| U_{p, 1}  (t, z )|  & \lesssim  B   e^{\la p \ra  \rho }  
 \int_0^z \!\!\!\! \int_0^t   e^{ -  p \rho'}    
 \frac{ k p \,    e^{ \rho'_p } }{\sqrt{(1 + \rho'_p)} }   dz' dt' 
 \\
 & \lesssim B  e^{\la p \ra  \rho }  
 \int_0^z \!\!\!\! \int_0^t       
   k p      e^{ - \mez p \rho' }     dz' dt'  
   \\
   & \lesssim
 B     \frac{1}{p | \uE |^2 }   \ln \big(1 +   \mez p^2 k |E|^2 zt \big)   
     \ e^{ \la p \ra \rho}  . 
 \end{aligned}
\end{equation*} 
On $\Omega $, 
$\mez k |E |^2  zt  \le   \rho^2  \le   | \ln \delta |^2 $. Thus, for $\delta \ge e^3 $,   
  Lemma~\ref{lem3.1b} implies that  
 \begin{equation}
\label{b4.14}
| U_{p, 1}  (t, z )|    \lesssim  \frac{1}{ | \uE |^2 }   \ln \big(1 +  | \ln \delta |^2 )    e^{ \la p \ra \rho}
    \lesssim     \frac{1}{ | \uE |^2 }  \ln ( |  \ln \delta |)  \   e^{ \la p \ra \rho}. 
\end{equation}
With \eqref{b4.10} and \eqref{b4.13}, this implies the estimate \eqref{b4.8}. 

\medbreak

{\bf b ) }   Similarly, when $p= 0$ the estimate for $n_0$ is immediate. 
When $p \ne 0$, by \eqref{2.27b}, $n_p$ is the sum of two terms. 
Up to constant factors, the first one is a convolution of $F_p$  by $\cE_1$
(computed at the frequency $kp$). Using again \eqref{a4.6} and \eqref{a4.8} 
this term satisfies 
\begin{equation}
\label{b4.15}
\begin{aligned}
| n_{p,1} (t, z )|  & \lesssim A     e^{p \rho (t, z) }  
 \int_0^z \!\!\!\! \int_0^t      
 \frac{ k p      e^{ - \mez p  \rho (t', z') } }{\sqrt{(1 + \rho'_p)} }   dz' dt' 
 \\
 & \lesssim  A    e^{\la p \ra  \rho }  
 \int_0^z \!\!\!\! \int_0^t       
   k p      e^{ - \mez p \rho' }     dz' dt' 
  \lesssim  A   \frac{1}{ | \uE |^2 }  \ln ( |  \ln \delta |) 
   e^{\la p\ra  \rho }    . 
 \end{aligned}
\end{equation} 
The second term in $n_p$ is the convolution of $h_p$ with $\cE_2$ and satisfies 
\begin{equation*}
| n_{p,2}  (t, z )|   \lesssim B   e^{\la p  \ra  \rho(t,z)  }  
 \int_0^z \!\!\!\! \int_0^t   e^{ -  p \rho'}    
 \frac{ k^2 p^2     e^{ \rho'_p } }{\tilde r'_p  \sqrt{(1 + \rho'_p)} }   dz' dt' 
\end{equation*}
with 
 $\frac{1}{\tilde r'_p}  \le   
  \big( t' + \frac{1}{kp} \big)$.  
 Thus, by Lemma~\ref{lem3.1} 
 \begin{equation}
\begin{aligned}
|n_{p,2}   (t, z )|  & \lesssim B    e^{\la p \ra  \rho }  
 \int_0^z \!\!\!\! \int_0^t       
   k^2 p^2     \big( \frac{1}{kp}   +  t ' \big)        e^{ - \mez p \rho(t', z') }    dz' dt' 
 \\
 &  \lesssim  B  \Big( \frac{k t}{| \uE |^2}   +    \frac{1}{ | \uE |^2 }  \ln ( |  \ln \delta |)    \Big)
   e^{p \rho }    . 
 \end{aligned}
\end{equation} 

\end{proof}



\section{Analytic  solutions of the Goursat problem }

Consider 
\begin{equation}
\label{x7.1}
\left\{ \begin{aligned}
& i  \D_z E + \Delta_x E = n E,\\
& (\D_t^2 - \Delta_x ) n = \Delta_x  \vert E \vert^2
\end{aligned}
\right.
\end{equation}
on $  t \in [0, T]$, $ z \in [0, Z] $ and  $x \in \TT^2$, together 
with initial boundary conditions 
\begin{equation}
\label{x7.2}
n _{| t = 0}  = 0, \quad \D_t n_{| t = 0} = 0, \quad  E_{ z = 0}  = E_0  . 
\end{equation}
   
   We prove that this Goursat problem is locally well posed in spaces 
   of analytic functions, following the general approach presented 
   in \cite{Wag}. One of our objective is to give an explicit 
   lower bound for the domain of existence, depending on the 
   norms of the boundary data (see Remark~\ref{rem56} below).

 \subsection{Spaces}
 
 Fix $ s > $ and equip $H^s(\TT)$ with a norm $\| \cdot \|_s$ such that 
 \begin{equation}
 \label{x7.3}
 \| u v \|_s \le \| u \|_s \ \|v\|_s. 
 \end{equation}
 \begin{defi}
 Given a (formal) power series $\phi ( \rx  ) = \sum \phi_n \rx^n$ in one variable 
 $\rx$, with 
 nonnegative coefficients, 
 we say that  $ u(x) \ll \phi $ if  for all $\alpha \in \ZZ^2$, 
 $\| \D_x^\alpha u  \|_s  \le  | \alpha | ! \,  \phi_{|\alpha|}$. 
 
 \end{defi} 
 
 \begin{lem}
 i)  If $u \ll \phi$, then $\D_{x_j} u \ll \phi'$. 
 
 \quad i) If $u \ll \phi$ and $v \ll \psi$, then $u v \ll \phi \psi $, 
 \end{lem}
 
 We will apply this definition to functions $u $ and power series which also depend 
 on $t$ and $z$, seen as parameters.  In particular, given a power series
 $\phi(\rx)$, we will consider power series
 \begin{eqnarray} 
\phi_{\lambda}  (    z, \rx  ) = \phi (  \rx  + \lambda z), 
 \end{eqnarray} 
 and use the notations $u(t,  z, x) \ll \phi_{ \lambda} (    z, \rx)  $ if 
 $u (t,   z, \cdot) \ll \phi_{  \lambda} ( z, \cdot)$.

 \begin{lem} 
    If $u (t, z, x) \ll    \phi'_{R, \lambda} (z, \rx)$ then 
\begin{equation*}
\int_0^z  u (t, z', x) dz' \ \ll \ \frac{1}{\lambda} \phi_{R, \lambda}(t,z,\rx). 
\end{equation*}
  
 \end{lem}

 In particular, we   consider the power series (see \cite{Wag})
 \begin{equation}
 \label{x7.5} 
 \varphi (\rx ) =  c_0 \sum_{ n= 0}^\infty \frac{ \rx^n}{ n^2 + 1} , 
 \quad   \varphi_{R, \lambda}  (z, \rx ) = \varphi (R \rx + \lambda z ). 
 \end{equation}
 where $c_0$ is such that $\varphi^2 \ll \varphi$.  
  We also introduce the notations 
 \begin{equation}
 \label{x.7.6}
  \varphi (R \rx + \lambda z ) = \sum_{n=0}^\infty \varphi_{n} (\lambda z)  R^n \rx^n. 
 \end{equation} 
 Of course the explicit value of $\varphi_n$ can be deduced from \eqref{x7.5}. 
 We also note that 
  \begin{equation}
 \label{x.7.7}
  \varphi'_{R, \lambda}( z, \rx)  = \sum_{n=0}^\infty (n+1)  \varphi_{n+1} (\lambda z)  R^{n+1} \rx^n. 
 \end{equation} 
 and by integration of $\varphi' ( R \rx + \lambda z)$ in $z$, we see that 
 \begin{equation}
 \int_0^z \varphi_{n+1} (\lambda z') dz'  = \frac{1}{ (n+1) \lambda} 
 \big(\varphi_n(\lambda z) - \varphi_n(0) \big)
 \le \frac{1}{ (n+1) \lambda} 
  \varphi_n(\lambda z)
 \end{equation}

 \begin{defi}
 Denote by $\EE_{R, \lambda} $ the space of functions $u$  on $[0, T] \times [0, Z] \times \TT^2$, with $Z = \frac{1}{\lambda} $, such that   
 \begin{equation}
 \label{x7.8}
 u (t, z, x) \ll  C \varphi_{R, \lambda}  (z, \rx)
 \end{equation}
 for some constant $C \ge 0$. 
 
 Similarly,   $\FF_{R, \lambda} $ denotes the space of functions $u$  on $[0, T] \times [0, Z] \times \TT^2$,   such that 
 \begin{equation}
 \label{x7.9}
 u (t, z, x) \ll  C \varphi'_{ R, \lambda} (z, \rx   )
 \end{equation}
 for some constant $C \ge 0$. 
 \end{defi}
 The norms in $\EE$ and $\FF$  are  the  best constant $C$ in \eqref{x7.8} and \eqref{x7.9}.

 For the boundary data $E_0 (t,x) $ on $[0, T] \times \TT^2$, we consider the spaces 
 $\EE_{R}$ of functions $u(t, x) \ll C \varphi( R \rx)$.

 \begin{theo}
 \label{thx7.5} 
 Let $E_0 \in \EE_{R} $. 
 Then, the problem \eqref{x7.1} has a unique solution 
 $(E, n) \in \EE_{R, \lambda} \times \FF_{R, \lambda} $ 
 if $\lambda \ge  C R T \| E_0 \|^2_{\EE_R}$ where 
 $C$ is a constant independent of $T$,  $R$ and $e_0$.

 \end{theo}

 \begin{rem}
 \label{rem56}
 \textup{Functions in $\EE_{R, \lambda} $  are defined  for 
 $R  | \im x | + \lambda | z |  < 1 $. Thus $ R^{-1}$  measures the width 
 of the complex domain where the boundary data is defined and $\lambda^{-1}$ 
 is the order of the length of propagation in $z$.  In particular, for boundary data
 $ E_0 (t,x)  =  e(t)  e^{ i k x} $  then, for $R \approx k$,   
 $\| E_0 \|_{\EE_R}  \lesssim   \uE   =  \| e \|_{L^\infty} $. 
 Theorem \ref{thx7.5} asserts that the length of stability 
 $Z = \lambda^{-1}$ satisfies }
 \begin{equation}
 \label{x7.9bb}
 Z T k  \uE^2  \lesssim 1 , 
 \end{equation}
 \textup{which is very similar to the condition $\rho \lesssim 1$ that was used in Section~3.  }
 
 \end{rem}

 
 \subsection{Resolution of \eqref{x7.1} in $\EE$}
 
 Consider the wave equation 
 \begin{equation}
 \label{x7.10}
 \D_t^2 n - \Delta_x n  = \Delta_x h , \quad  n_{| t = 0} = \D_t n_{| t = 0}  = 0. 
 \end{equation}
 
 \begin{lem}
 For $h \in \EE_{R, \lambda}$ , the solution $n $ of \eqref{x7.10} belongs to
 $\FF_{R, \lambda}$  and 
 \begin{equation}
 \| n \|_{\FF_{R, \lambda}}   \le  C_0  T     \| h \|_{\EE_{R, \lambda}} 
 \end{equation} 
 where $C_0$ is independent of $T$, $R$ and $\lambda$. 
 \end{lem}
 
 \begin{proof}
 The mean value of $n$ (the $0$-th Fourier coefficient) vanishes and therefore
 the $H^s$ norm of $n$ is equivalent to the $H^{s-1}$ norm of $\D_x n$. 
 Thus, by standard energy estimates for the wave equation, there holds
 \begin{equation*}
 \| \D_x^\alpha n (t, z, \cdot )   \|_s  \le C \int_0^t \| \D_x^\alpha h(t', z, \cdot )  \|_{s+1} dt' . 
 \end{equation*}
 Thus, for $| \alpha | = n$, 
  \begin{equation*}
 \| \D_x^\alpha n (t, z, \cdot )   \|_s  \le C T \| h \|_{\EE_{R, \lambda}} \  (n+1) !   \ 
  R^{n+1}  \varphi_{n+1} (\lambda z)  .  
 \end{equation*}
  Therefore
  \begin{equation*}
  n (t, z, x) \ll   C T \| h \|_{\EE_{R, \lambda}} \  \varphi'_{ R, \lambda} (z, \rx   )
  \end{equation*}
 and the lemma follows. 
 \end{proof}
 
 We consider next the Schr\"odinger equation 
 \begin{equation}
 \label{x7.11}
 i \D_t e  - \Delta_x e  = f  , \quad  e_{| z = 0} = e_0. 
 \end{equation}
 
  \begin{lem}
 For $f \in \FF_{R, \lambda}$, and $ e_0 \in \EE_{R}$, the solution $e $ of \eqref{x7.11} belongs to
 $\EE_{R, \lambda}$  and 
 \begin{equation}
 \| e \|_{\EE_{R, \lambda}}   \le  \frac{  R }{\lambda}     \| f \|_{\FF_{R, \lambda}}  + \| e_0 \|_{\EE_R} .
 \end{equation} 
 \end{lem}
 \begin{proof}
 Standard energy estimates  imply that 
  \begin{equation*}
 \| \D_x^\alpha e (t, z, \cdot )   \|_s  \le   \int_0^z \| \D_x^\alpha  f (t, z', \cdot )  \|_{s} dz' 
 + \| \D_x^\alpha e_0 (t, \cdot )   \|_s
 \end{equation*}
 Thus, for $| \alpha | = n$, 
  \begin{equation*}
 \begin{aligned}
 \| \D_x^\alpha e (t, z, \cdot )   \|_s  \le   \| h \|_{\FF_{R, \lambda}}  \  (n+1) !  
 R^{ n+1}    \ 
   \int_{0}^z  & \varphi_{n+1} (\lambda z') dz'  
  \\ =
  & \| h \|_{\FF_{R, \lambda}} \  n ! \   \frac{R^{ n+1}}{\lambda}   \ 
    \varphi_{n} (\lambda z)  . 
   \end{aligned} 
 \end{equation*}
  Therefore
  \begin{equation*}
  e (t, z, x) \ll   \frac{R}{\lambda}  \| f \|_{\FF_{R, \lambda}} \  \varphi_{ R, \lambda} (z, \rx   )
  \end{equation*}
 and the lemma follows. 
 
 \end{proof}

 \begin{lem}
 For $e_1 \in \EE_{R, \lambda}$, $e_2 \in \EE_{R, \lambda}$ and $n \in \FF_{R, \lambda}$, there holds
 $ e_1 e_2  \in \EE_{R, \lambda}$ and   $ n e_1 \in \FF_{R, \lambda}$ with 
 \begin{eqnarray*}
 && 
 \| e_1 e_2 \|_{\EE_{R, \lambda}} \le   \|  e_1   \|_{\EE_{R, \lambda}}
 \|  e_2   \|_{\EE_{R, \lambda}}, 
 \\
 &&
  \| n e _1 \|_{\FF_{R, \lambda}} \le   \|  e_1   \|_{\EE_{R, \lambda}} \   \| n  \|_{\FF_{R, \lambda}} . 
 \end{eqnarray*}

 \end{lem} 
 
 \subsection{Proof of Theorem~\ref{thx7.5}}
 
 Denote by $\cN(h)$ the solution of \eqref{x7.10} and by $\cS_0 (e_0) + \cS (f)$  the solution 
 of \eqref{x7.11}. To prove Theorem~\ref{thx7.5} it is sufficient to show that there is $C > 0$, such that  for  
$ \lambda \ge C R T \| E_0 \|^2_{\EE_R} $, the equation 
\begin{equation}
\label{x7.14}  E = \cS(E_0)  + \cT (E, E, E)  
\end{equation}
 has a unique solution in $\EE_{R, \lambda}$, where $\cT$ is the trilinear operator 
 \begin{equation*}
 \cT ( u, v, w) = \cS \big( u \cN(v\overline w) \big) . 
 \end{equation*}
 The estimates above show that $\cT$ maps 
 $(\EE_{R, \lambda} )^3$ to $\EE_{R, \lambda}$ and that
 \begin{equation*}
  \| \cT (u, v, w)  \|_{\EE_{R, \lambda}} \le  \frac{C_0 R T}{\lambda} 
   \| u \|_{\EE_{R, \lambda}}  \| v \|_{\EE_{R, \lambda}}  \| w \|_{\EE_{R, \lambda}} . 
 \end{equation*}
From here, standard Picard's iterates imply Theorem~\ref{thx7.5}.


\section{Remarks and comments} 
 
 \subsection{Nonconstant backgrounds} 
 
 The boundary data we have considered in Sections~2 and 3  are of the form
 \begin{equation}
 \label{6.1.0}
 1_{ \{ t \ge 0 \}}\, \big(  \uE  +    e (t, x) \big), 
 \end{equation} 
 thus have a jump at $t = 0$, which is not physical.  Note that arbitrary functions 
 $\uE(t)$ together with $n = 0$ are solutions of  \eqref{1.3}. One could take them 
 as background state. We now briefly sketch how one can start the analysis. 
  Again, there is no restriction  in assuming that $\uE$ is real.  
  The linearized systems \eqref{2.3} and \eqref{2.4} are  unchanged, except that 
  $\uE$ is now a function of time.   Therefore, we now proceed by using 
the   Fourier-Laplace in $z$ rather than in time:  we extend $u, v, n$ by $0$ for 
 $z < 0$. Note that the extension $\tilde u$ satisfies 
 $\D_z \tilde u = \widetilde {\D_z u} + u_0 \delta_{ z = 0}$. 
 Therefore, the Fourier-Laplace transforms satisfy 
   \begin{equation}
\label{6.1.1}
\left\{ \begin{aligned}
& - (\zeta + k^2  )  \hat u  -    \uE   \hat n =  i u_0    ,\\
& ( - \zeta     +  k^2)   \hat   v   +     \uE   \hat n =  i v_0     , \\ 
& (\D_t^2 + k^2 )  \hat n + k^2    \uE   (\hat u  + \hat v  ) = 0  . 
\end{aligned}
\right.
\end{equation}
Hence, 
   \begin{equation}
\label{6.1.2}
\left\{ \begin{aligned}
&   \hat u   = -    \frac{ \uE   \hat n}{ \zeta + k^2}  - \frac{  i u_0} { \zeta + k^2}     ,\\
&   \hat   v   =      \frac{ \uE   \hat n}{ \zeta - k^2}  -   \frac{  i v_0 }{ \zeta - k^2}     , \\ 
\end{aligned}
\right.
\end{equation}
and 
\begin{equation}
\label{6.1.3}
  \big(\D_t^2 + k^2   (1 +  a p(t)    \big) \hat n  =   k^2  \hat h , 
\end{equation}
 with 
 \begin{equation}
 \label{6.1.4}
 a  =    |E |^2 \ \frac{ 2 k^2  }{ \zeta^2 - k^4} , \qquad p(t) = | \uE (t)|^2 
 \end{equation}
 \begin{equation}
 \label{6.1.5}
 h  =   i E  \big( \frac{ u_0   }{ \zeta - k^2}   + \frac{ v_0 }{ \zeta + k^2} \big). 
 \end{equation}
 
 Denote by $\cN (t, s, a)$ the fundamental solution of the second order o.d.e. 
 \begin{equation}
 \big( \D_t^2 + k^2 (1 + a p) \big) \cN  = 0, \qquad \cN _{| t = s} =0, \quad 
  \D_t \cN _{| t = s} = 1. 
 \end{equation}
  Thus 
 \begin{equation}
 \hat n(t) =  k^2 \int_0^t  U(t, s, a)   h(s)  ds. 
 \end{equation}
 Performing the inverse Laplace transform, we get
 \begin{equation}
 n  =    \int_0^t \cE_1 (t, s, z)  u_0 (s) ds   + \int_0^t \cE_2 (t, s, z)  v_0 (s) ds 
 \end{equation} 
 with 
 \begin{equation}
 \label{6.1.10}
  \cE_1(t, s, z)  =  \frac{i \uE(s) }{2\pi} \int_{\RR - i \gamma} e^{ i z \zeta} U (t, s, a ) \frac{ k^2 \, d \zeta}{ \zeta - k^2} , 
 \end{equation}
  \begin{equation}
  \label{6.1.11}
  \cE_2(t, s, z)  = \frac{i \uE(s)}{2\pi}  \int_{\RR - i \gamma} e^{ i z \zeta} U (t, s, a ) \frac{ k^2 \, d \zeta}{ \zeta + k^2} . 
 \end{equation}

 \begin{lem}
 \label{lemxxx2.1}
 $\cN$ is an entire function of  $a$. 
 Moreover, for   $| a | $small and  $0 \le s \le t \le 1$ there holds  
 \begin{equation}
 \cN(t,s, a)    =   \frac{1}{k} \Big(  \alpha (t, s,a) e ^{ i k \vp (t, s, a)} - \beta (t,s,a ) e^{ - i k \vp(t, s, a)} \Big) 
 \end{equation}
with  
 \begin{equation}
 \D_t \vp =    \sqrt{ 1 + a p}  , \quad \vp_{| t = s }  = 0, 
 \end{equation}
 and $\alpha$ and $\beta$ are $O(1)$ as well as their time derivatives. 
 \end{lem}

 Because $| a | \to 0$ as $|\zeta | \to \infty$, $U$ is holomorphic in $\zeta$ and bounded 
 for large $\zeta$. Thus, one can close the integration path in \eqref{6.1.10} and \eqref{6.1.11}: 
 for instance,   for $z > 0$: 
 \begin{equation}
  \cE_1(t, s, z)  =   \frac{1}{2\pi}\int_{\Gamma} e^{ i z \zeta} U (t, s, a ) \frac{ k^2 \, d \zeta}{ \zeta - k^2} 
 \end{equation}
 where $\Gamma$ is a large circle  in the complex plane. 
  
 From here, using the Lemma above, one can repeat most of the computations performed in Section~2, 
 with only technical additional difficulties. The important point is that the 
 amplification is now driven by the kernel $e^\rho$ with
  \begin{equation}
  \rho := \sqrt{ 2 k (z- z')  (P(t) - P(t') )} , \qquad  P(t) = \int_0^t | \uE(s) |^2 ds . 
  \end{equation}  
  Of course, when $\uE$ is constant we recover the previous results. 
  On an interval where $ |\uE | > 0$,  there holds 
  $  P(t) - P(t') \approx  (t- t')$ and all the estimates of Section 2 are unchanged. 
 
 In this direction,  another interpretation of the constant case $\uE$,  different from 
 \eqref{6.1.0} is that the boundary data is 
  \begin{equation}
 \label{6.2.1}
   \uE  (t)   +      e (t, x)   1_{ \{ t \ge 0 \}}\,  
 \end{equation} 
 where $\uE $ vanishes for $t \le -T_0 < 0$ and is equal to the constant $\uE$ for $t \ge 0$. 
 This means that the perturbation $e$ starts when the background is stabilized at a constant value. 
 This is  not physically realistic either,  but the expectation is that   this analysis is qualitatively relevant for physical interpretations.

 
 \subsection{The amplification rate} 
 
 We come back to the case where $\uE$ is constant.  
 The solution of the homogeneous linearized equations \eqref{2.4} is given in Lemma~\ref{lem2.1}. 
 It involves convolution by functions defined by contour integrals \eqref{2.14}. 
 These integrals were estimated in Section~2 using deformations of contours and the saddle point 
 method. According to Lemma~\ref{lem2.4f}, 
 a model for  the integral $ \cE_1$,  is 
  \begin{equation}
 \label{6.2.2} \cF_1   = \frac{1}{2i  \pi} \int_\Gamma e^{ i ( t \zeta - z k \lambda) } \frac{ d \zeta} { \zeta - k } , \qquad \lambda =  k  +  \frac{ | \uE |^2}{ 2 (\zeta - k)} . 
 \end{equation}
 Note that 
 \begin{equation}
 \label{6.2.3}
 \cF_1 = e^{ i ( t k - z k^2)}  J (\rho), \qquad \rho = \sqrt{ 2 k  | \uE |^2 z t }
 \end{equation}
 where $J$ is the bessel function
 \begin{equation*}
 J (\rho )  =  \frac{1}{2i  \pi} \int_\Gamma e^{ i \mez  \rho ( \zeta - \frac{1}{\zeta} ) } \frac{ d \zeta} {\zeta}
 = \int_0^{2\pi}  e^{ \rho \sin \theta}  d \theta  = \sum_{n=0}^\infty \frac{\rho^n}{n!^ 2}. 
 \end{equation*}
 This is indeed the core of the analysis performed in Section~2.  
 
 Following this path, a model for the solutions of the homogeneous equation \eqref{2.4} is
 \begin{equation}
 \label{6.2.4}
 u (t, z) = \int_0^t e^{ i ( (t- t')  k - z k^2)}  J (\sqrt{ 2 k  | \uE |^2 z (t-t') } )  u_0 (t') dt' .    
 \end{equation}
 The kernel has an amplification part given by $J$ but it also has an oscillating part 
 given by $e^{ i  (t- t')  k}$.  In Theorem~\ref{th2.7}, we have chosen $u_0$ such that 
 it has itself an oscillation $e^{ i k t'}$ to eliminate the oscillations in \eqref{6.2.4}  and maximize 
 the amplification. 
 For general $u_0$ the actual amplification results from a delicate balance 
 between $J$ and the oscillations. Typically, one integration by part in the integral 
 above, wins a factor $k^{-\mez}$ and $k^{ -1} $ times a derivative of $u_0$. 
 Thus the amplification $e^{\rho}$  is cut down by a factor $k^{- m}$  where
 $m$ is the smoothness of $u_0$. 
 This can also be seen on the Fourier Laplace representation of $u$: 
 \begin{equation}
 \label{6.2.5}
 u(t, z)    = \frac{e^{ i ( t k - z k^2)} }{2i  \pi} \int_\Gamma e^{ i ( t \zeta - \frac{z k|\uE |^2}{\zeta}   ) } 
 \hat u_0  (k + \zeta)  \frac{  d \zeta} { \zeta  } .  
 \end{equation}
In any case,  there is one clear  conclusion : \emph{the more amplified $u_0 $ are those who have a nontrivial 
 oscillation at frequency $k$}.

 \bigbreak
 
 This can also be seen in a different way. Consider the Fourier-Laplace 
 expression of the solution: 
 \begin{equation}
 \label{6.2.6}
 \hat U (z, \zeta)  = e^{ i z A (\zeta) } \hat U_0 (\zeta) . 
 \end{equation}
For $\zeta = ( \tau - i \gamma)  \notin \{  -k, + k \}$,  the eigenvalues of 
 $A$ are 
 \begin{equation}
 \label{6.2.7}
\mu (\zeta) =  \pm  k^2  \sqrt { 1 + \frac{1}{b}  } , \qquad  b = \frac{ \zeta^2 - k^2}{ | \uE |^2}. 
 \end{equation} 
 Consider    \emph{real} time frequencies 
 $\zeta = \tau \notin \{  -k, + k \}$.   
 When $1 + \frac{1}{b} > 0$, the eigenvalues are real and there is no amplification. 
  When $1 + \frac{1}{b} < 0$, that is when 
 \begin{equation}
 \label{6.2.8}  0 <  k^2 - \tau^2 < | \uE |^2 , 
 \end{equation}
 the eigenvalue is purely imaginary and this suggests that the  \emph{amplification factor for 
 oscillations with frequency $\tau$} is 
 \begin{equation}
 \label{6.2.9} 
 e^{ z  | \im \mu |   }     , \qquad  |  \im \mu |  =  k^2 \frac{\sqrt{| \uE |^2 - k^2 + \tau^2}}{\sqrt{  k^2 -  \tau^2}} . 
 \end{equation}
 In particular,  only real frequencies close to $k$ are amplified.

 Of course, considering only real frequencies is not consistent with the forward evolution problem
 under consideration and the condition that $U_0$ vanishes in the past. 
 $e^{ i z A (\zeta)}$ has essential singularities at $+ k$ and $-k$, and one has to consider
 integration on complex path as in Section~2, to turn around the singular points.  
 Then the imaginary part of the phase to consider is 
 $t \im \zeta - z \im \mu $. In Section~2 we have chosen  optimal  contours in order to obtain   sharp  estimates.  In particular, the amplification \eqref{6.2.4} is not correct. But the computation suggests
  that  the amplification is mainly due to frequencies close to $k$, and this is correct.


\subsection{The physical  context  and physical values}

The  Zakharov's
equations \cite{zakharov} have been introduced  at the beginning of the 70's,  
to describe  electronic plasma
waves. They  couple the slowly varying envelope of the electric 
field  and  the low-frequency variation of
the density of the ions.  When modeling the propagation of a laser beam  in a plasma,  
several phenomena occur. One has   to take into account   the laser beam itself,   the Raman component and  the electronic plasma waves
(see \cite{colins1,colins2} for example). 
For laser propagation or for the Raman component, one often  uses the paraxial 
approximation and the Zakharov system that couples the envelop of the   
vector potential $A$ of the electromagnetic field $A$  to the low-frequency variation of
the density $\delta n$ of the ions reads
\begin{equation}
\label{x8.1}
\left\{ \begin{aligned}
& i( \D_t  + \frac{k_0 c^2}{\omega_0}\D_z)A + \frac{c^2}{2\omega_0}\Delta_x A =\frac{\omega_{pe}^2}{2n_0\omega_0} \delta n\,  A,\\
& (\D_t^2 - c_s^2\Delta_x) \delta n = \frac{\omega_{pe}^2}{4\pi m_i c^2}\Delta_x \vert A \vert^2,
\end{aligned}
\right.
\end{equation}
where $\omega_0$ is the frequency of the laser, $k_0$ its wave number and $\omega_{pe}$ the plasma electronic frequency; they  are linked by
the dispersion relation  $\omega_0^2 = \omega_{pe}^2+k_0^2c^2$ where   $c$ is the speed of light;  $n_0$ the mean density of the plasma, $m_i$ is the mass of the ions and $c_s$ the sound
velocity in the plasma.  In suitable units for $n_0$, the plasma frequency $\omega_{pe}$ is given 
by
\begin{equation}
\label{x8.2a}
\omega_{pe} = \sqrt{ \frac{4 \pi e^2 n_0}{m_e}} 
\end{equation}
where $- e$ is the charge of the electron  and $m_e$ its mass. 

As above,  $z$ is the space component in the direction of propagation of the laser beam
 and $x$ denotes  the space components in directions that are transversal  to the propagation.  
 
Introducing the dimensionless quantities   
\begin{equation}
\label{x8.2.b}
 n = \delta n / n_0, \qquad   a =  \frac{A}{\underline A} , 
 \quad \underline A = \frac{m_e  c c_s}{e} 
 \end{equation}
  and using \eqref{x8.2a}, 
 the system reads
\begin{equation}
\label{x8.3a}
\left\{ \begin{aligned}
& i(  \frac{\omega_0}{c^2} \D_t  +  k_0 \D_z) a  + \mez \Delta_x a =\frac{\omega_{pe}^2}{2 c^2 } \,  n\,  a,\\
& ( \frac{1}{c_s^2} \D_t^2 -   \Delta_x)  n = \frac{m_e}{m_i} 
 \Delta_x \vert a  \vert^2,
\end{aligned}
\right.
\end{equation}

In addition to the pulsation $\omega_0$ and the wave length 
$\lambda_0 = \frac{2\pi}{k_0}$, characteristic quantities for the laser beam are 
its duration $\tau_0$, its transversal width $R_0$ (supposed to be smaller than the transversal dimension of the plasma), the expected length of propagation $Z_0$  
and also the characteristic dimension of the transversal variations (speckles)
$\lambda$ which correspond to variations of  the spectrum by frequencies 
$k \approx  \frac{2\pi}{\lambda}$.  In the application we have in mind
$\lambda \ll R_0$  and  both for theoretical and computational  reasons (spectral methods)  
it makes sense to \emph{assume periodicity in $x$} with period $X$ satisfying
$\lambda \le X \ll R_0$. 

With this data in mind, rescale  the space-time variables introducing the characteristic transversal 
width $X$: 
\begin{equation}
\label{x8.4}
\tilde x = \frac{ x}{X }, \quad 
 \tilde z = \frac{ z}{2 k_0  X ^2}, \quad \tilde  t = \frac{  c_s t }{X }
\end{equation}
In these variables, the system reads
\begin{equation}
\label{x8.6}
\left\{ \begin{aligned}
& i( \epsilon  \D_{\tilde t}  +    \D_{\tilde z} ) a +    \Delta_{\tilde x} a =
\alpha ^2   n  a ,\\
& ( \D_{\tilde t}^2 -  \Delta_{\tilde x}) n =  \frac{m_e}{m_i}  \Delta_{\tilde x} \vert a \vert^2,
\end{aligned}
\right.
\end{equation}
with 
\begin{equation}
\epsilon = \frac{ 2  \omega_0 c_s X   }{c^2}, \qquad 
\alpha =  \frac{\omega_{pe}X }{  c}. 
\end{equation} 
With 
\begin{equation}
\label{x8.7}
\tilde n =   \alpha^2  n , \qquad  \tilde a  =  \alpha \sqrt{ \frac{m_e}{m_i}}  a 
\end{equation}
we obtain the dimensionless system: 
\begin{equation}
\label{x8.8}
\left\{ \begin{aligned}
& i( \epsilon  \D_{\tilde t}  +    \D_{\tilde z} ) \tilde  a +  \Delta_{\tilde x} \tilde a =
\tilde  n \tilde a,\\
& ( \D_{\tilde t}^2 -  \Delta_{\tilde x}) \tilde n =    \Delta_{\tilde x} \vert \tilde a \vert^2. 
\end{aligned}
\right.
\end{equation}

For this system, the typical \emph{amplification factor}  is 
\begin{equation}
\label{x8.9}
  \rho   =  \sqrt { 2 | \tilde a |^2 \, \tilde k\,  \tilde \tau_0 \,   \tilde Z_0 }  
\end{equation}
where $\tilde k$ is the scaled frequency of speckles,  
$\tilde \tau_0 $ the scaled  time of propagation under consideration and 
$\tilde Z_0$ the scaled length of propagation. 
Scaling back to the original variables, there holds 
\begin{equation}
\tilde Z_0  = \frac{ Z_0}{2 k_0 X^2  } , \quad  \tilde k = X  k, \quad 
 \tilde \tau = \frac{ c_s \tau_0 }{X}, 
\end{equation}
thus

\begin{prop}  The dimensioned amplification factor $\rho$ of Section~$2$ is given by 
\begin{equation}
\label{x8.11} 
  \rho^2  =   2  \frac{|   A |^2}{  \underline A  ^2}     \frac{m_e  c_s \omega^2_{pe} }{   m_i    c^2 }  \ \frac{ k }{  k _0 } \, \tau_0 Z_0 . 
\end{equation}
\end{prop}

\bigbreak
For the physical significance of the analysis it is important to 
evaluate the various constants. 
The velocity of light is $c = 3\   10^8 \mathrm{ m s^{-1}}$, the sound velocity of the electrons 
is of order $c_s \approx  0.005 c = 1.5 \ 10^6 \mathrm{ m s^{-1}} $. 
With $ e = 1.6 \ 10^{-19} \mathrm{C}$  and $m_e = 0.9\  10^{-30} \mathrm{Kg}$, 
the scaling factor for $A$ is 
\begin{equation}
\underline A = \frac{m_e  c c_s}{e}   \approx 2.5 \, 10^{3} \, \mathrm{V}  . 
\end{equation}
Typical values of  
$\omega_{pe}$ are of order  $ 10^{15} s^{-1}$ and $\frac{m_e}{m_i} \approx   10^{-4}$.  
Thus 
\begin{equation}
2  \frac{m_e  c_s \omega^2_{pe} }{   m_i    c^2 } \approx  4  \   10^{15}  \     \ \mathrm{ m ^{-1}s^{-1} }
\end{equation} 
A typical value of  $\lambda_0$ is  $\lambda_0 = 0.35 \mu\mathrm{m}$
corresponding to a wave number 
$k_0 = \frac{2 \pi}{ \lambda_0} =  1.8 \  10^7 \mathrm{m^{-1}}$. The pulsation is given 
by the dispersion relation   $\omega_0^2 = \omega_{pe}^2+k_0^2c^2$, 
yielding 
$\omega_0 =  5.5 \ 10^{15}  \mathrm{ s^{-1}}$.

By construction, the relative variation $n = \delta n / n_0$ is small compared to $1$
and the consistency of  \eqref{x8.3a} requires that 
$ \frac{m_e}{m_i} | a |^2$ must also be small compared to $1$, that is 
$| a | \ll 100 $ yielding the upper bound for $A$
\begin{equation}
\label{x8.26} 
| A | \ll  10^5 \ \mathrm{V}. 
\end{equation} 
  The envelop  of the electromagnetic field is linked to $A$ by  the polarization relation
  \begin{equation}
  E  =  i \frac{ \omega_0}{c}  A  
  \end{equation}
  yielding the bound
  \begin{equation}
  | E | \ll  10^{12}  \ \mathrm{V m^{-1}},
  \end{equation} 
which allows for very intense fields.

 Typical values for the dimension of the speckles is $\lambda \approx 10  \mu \mathrm{m}$, 
 corresponding to frequencies $k = \frac{2\pi}{\lambda} \approx 6 \ 10^{5} \mathrm{m^{-1}} $.  A constraint for the validity of the paraxial approximation which sustains the 
 Schr\"odinger envelop equation in \eqref{x8.1} is that 
 \begin{equation}
 k  Z_0  \le k_0  R_0 
 \end{equation}
 With $ k_0 / k \approx  30 $   as above  and  a diameter $R_0 $ ranging from $ 1 \mathrm{mm}$
 to  $ 1 \mathrm{cm}$, this allows for 
 propagations along distances $Z_0$ of order  $1$ to several $ \mathrm{cm}$. 
In this range of $k$, we finally obtain: 
\begin{equation}
\label{x8.29.0}
\rho =  \frac{ |A |}{ \underline A}  \sqrt {   \beta    \tau_0  Z_0 } \qquad 
\mathrm{with} \ \ \beta   \approx     10^{14}  \     \ \mathrm{ m ^{-1}s^{-1} }
\end{equation}

 The duration of the laser impulsion 
 $\tau_0 $  may be of order $10^{- 8} $ to $ 10^{-10} \mathrm{s}$.  Note that 
 for a $ 1 \mathrm{cm}$  long propagation of the laser beam, the time elapsed is of order  
 $3  10^{-11} \mathrm{s}$. 
 For $Z_0 $ of order  $1 $cm  and  $\tau_0 $ ranging from $10^{- 8} $ to $ 10^{-10} \mathrm{s}$, we  obtain that 
   \begin{equation}
\label{x8.29.1}
\rho =  \gamma  \frac{ |A |}{ \underline A}   
\end{equation}
 with $\gamma $ ranging from $10 $ to $ 100 $. 
 
 Thus, for $\rho$ to  be small, one must have
 \begin{equation}
 \label{6.3.45}
 | A |   \lesssim 10 \, \mathrm{V}, \quad \mathrm{or} \quad | E | \lesssim     10^8 \, \mathrm{V m^{-1}}, 
 \end{equation} 
 which does nor correspond to high intensity beams.

 Note also that with $X \approx 50 \ \mu\mathrm{m}$, there holds
\begin{equation}
 \epsilon \approx   9,  
\qquad \tilde \tau \approx 3. 
\end{equation}

\bigbreak
In conclusion, we see that  the values obtained in \eqref{6.3.45}  are compatible with the 
paraxial approximation and the derivation of the model \eqref{x8.1}.  
However, they do not allow for high intensity or high energy beams. 
This suggests that the model could be unadapted to such situations.


\vfill \eject


\end{document}